\documentclass[hidelinks,onefignum,onetabnum]{siamart220329}

\usepackage{lipsum}
\usepackage{amsfonts}
\usepackage{graphicx}
\usepackage{epstopdf}
\usepackage{algorithmic}
\ifpdf
  \DeclareGraphicsExtensions{.eps,.pdf,.png,.jpg}
\else
  \DeclareGraphicsExtensions{.eps}
\fi

\newsiamremark{remark}{Remark}
\newsiamremark{hypothesis}{Hypothesis}
\crefname{hypothesis}{Hypothesis}{Hypotheses}
\newsiamthm{claim}{Claim}

\headers{deep density approximation for stochastic systems}{Junjie He, Qifeng Liao, and Xiaoliang Wan}

\title{Adaptive deep density approximation for stochastic dynamical systems
\thanks{Submitted to the editors DATE.
\funding{
The first two authors are supported by the National Natural Science Foundation of China (No. 12071291). The third author is supported by NSF grant DMS1913163.
}
}
}

\author{Junjie He\thanks{School of Information Science and Technology, ShanghaiTech University, Shanghai, Shanghai 201210, China 
  (\email{hejj1@shanghaitech.edu.cn}).} 
  \and Qifeng Liao\thanks{Corresponding author. School of Information Science and Technology, ShanghaiTech University, Shanghai, Shanghai 201210, China (\email{liaoqf@shanghaitech.edu.cn}).} 
  \and Xiaoliang Wan\thanks{Department of Mathematics and Center for Computation and Technology, Louisiana State University, Baton Rouge, LA 70803 USA (\email{xlwan@math.lsu.edu})}.}

\usepackage{amsopn}
\DeclareMathOperator{\diag}{diag}

\ifpdf
\hypersetup{
  pdftitle={Adaptive deep density approximation for stochastic dynamical systems},
  pdfauthor={Junjie He, Qifeng Liao, and Xiaoliang Wan}
}
\fi

\usepackage{amsmath,amsfonts,amssymb}
\usepackage{subfigure}

\DeclareMathOperator{\Var}{Var}
\newcommand{\SiLU}{\operatorname{SiLU}}
\newcommand{\subT}[1]{\mathbf{T}^{(#1)}}
\newcommand{\resx}{\boldsymbol{x}_{\text{res}}}
\newcommand{\rest}{{t}_{\text{res}}}
\newcommand{\txtval}{\text{val}}
\newcommand{\bx}{\boldsymbol{x}}

\newcommand{\by}{\boldsymbol{y}}
\newcommand{\bz}{\boldsymbol{z}}
\newcommand{\bh}{\boldsymbol{h}}
\newcommand{\ba}{\boldsymbol{a}}
\newcommand{\bb}{\boldsymbol{b}}
\newcommand{\bX}{\boldsymbol{X}}
\newcommand{\bZ}{\boldsymbol{Z}}

\newcommand{\mathf}{\boldsymbol{f}}
\newcommand{\mathg}{\boldsymbol{g}}
\newcommand{\maths}{\boldsymbol{s}}
\newcommand{\matht}{\boldsymbol{t}}

\newcommand{\bxi}{\boldsymbol{\xi}}

\newcommand{\dif}{\mathrm{d}}

\newcommand{\bW}{\mathbf{W}}

\newcommand{\bbeta}{\boldsymbol{\beta}}
\newcommand{\bphi}{\boldsymbol{\phi}}
\newcommand{\mathN}{\mathbf{N}}
\newcommand{\mathT}{\mathbf{T}}

\newcommand{\mathI}{\mathbb{I}}

\newcommand{\adaptive}{\text{adaptive}}
\newcommand{\txttime}{\text{time}}

\newcommand{\txtinterface}{\text{interface}}

\renewcommand{\algorithmicrequire}{ \textbf{Input:}}     
\renewcommand{\algorithmicensure}{ \textbf{Output:}}    
\DeclareMathOperator*{\argmin}{arg\,min}

\usepackage{xcolor,color}
\usepackage{comment}

\begin{document}

\maketitle

\begin{abstract}
In this paper we consider adaptive deep neural network approximation for stochastic dynamical systems. Based on the Liouville equation associated with the stochastic dynamical systems, a new temporal KRnet (tKRnet) is proposed to approximate the probability density functions (PDFs) of the state variables. The tKRnet gives an explicit density model for the solution of the Liouville equation, which alleviates the curse of dimensionality issue that limits the application of traditional grid based numerical methods. To efficiently train the tKRnet, an adaptive procedure is developed to generate collocation points for the corresponding residual loss function, where samples are generated iteratively using the approximate density function at each iteration. A temporal decomposition technique is also employed to improve the long-time integration. Theoretical analysis of our proposed method is provided, and numerical examples are presented to demonstrate its performance. 
\end{abstract}



\section{Introduction}
\label{sec:intro}
Stochastic dynamical systems naturally emerge in simulations and experiments involving complex systems (e.g., turbulence, semiconductors, and tumor cell growth), where probability density functions (PDFs) of 
their states are typically governed by partial differential equations (PDEs)  \cite{lord2014introduction,cho2016numerical,tartakovsky2017method,brennan2018data}. These include the Fokker-Planck equation \cite{moss1989noise,risken1996fokker}, which allows assessing the time evolution of  PDFs in Langevin-type stochastic dynamical systems driven by Gaussian white noise, and  Liouville equations \cite{sobczyk1991stochastic,klyatskin2005dynamics,li2009stochastic,villani2009Optimal} which model the evolution of PDFs in stochastic systems subject to random initial states and input parameters. However, it is challenging to solve these PDEs efficiently due to difficulties such as high-dimensionality of state variables, possible low regularities, conservation properties, and long-time integration \cite{cho2016numerical}. This paper is devoted to developing new efficient deep learning methods for the Liouville equations to address these issues.

The main idea of deep learning methods for solving PDEs is to reformulate a PDE problem as an optimization problem 
and train deep neural networks to approximate the solution by minimizing the corresponding loss function. Based on this idea, many techniques have been investigated to alleviate the difficulties existing in applying traditional grid methods (e.g., the finite element methods \cite{elman2014finite}) for complex PDEs. These include, for example, 
deep Ritz methods \cite{weinan2017proposal,e2018deep}, physics-informed neural networks (PINNs) \cite{raissi2019physics, karniadakis2021physics}, deep Galerkin methods \cite{sirignano2018deep}, Bayesian deep convolutional encoder-decoder networks \cite{zhu2018bayesian,zhu2019physics}, weak adversarial networks \cite{zang2020weak} and deep multiscale model learning \cite{wang2020deep}. Deep neural network methods for complex geometries and interface problems are proposed in \cite{sheng2021pfnn, gao2021phygeonet, wang2020mesh}, and domain decomposition based deep learning methods \cite{li2020deep,jagtap2020conservative,li2020deep1, heinlein2020combining, kharazmi2021hp, dong2021local,xu2024domain} are studied to further improve the computational efficiency. In addition, residual network based learning strategies for unknown dynamical systems are presented in \cite{wu2020data,chen2022deep}. 

As the solution of the Liouville equation is a time-dependent probability density function, solving this problem can be considered as a time-dependent density estimation problem. While density estimation is a central topic in unsupervised learning \cite{scott2015multivariate}, we focus on the normalizing flows \cite{dinh2017density,kingma2018glow,chen2018neural,Kobyzev2021flowreview} in this work. The idea of the normalizing flows is to construct an invertible mapping from a given simple distribution to the unknown distribution, such that the PDF of the unknown distribution can be obtained by the change of variables. Constructing an efficient invertible mapping is then crucial for normalizing flows. In our recent work \cite{tang2020deep},  based on the Knothe-Rosenblatt (KR) rearrangement \cite{carlier2010KnotheTransport} and a modification of affine coupling layers in real NVP \cite{dinh2017density}, a normalizing flow model called KRnet is proposed. A systematic procedure to train the KRnet for solving the steady state Fokker-Planck equation is studied in \cite{TANG2022111080}. In addition, an adaptive learning approach based on temporal normalizing flows is proposed for solving time-dependent Fokker-Planck equations in \cite{feng2021solving}.

In order to efficiently solve the Liouville  equation, we generalize our KRnet to time-dependent problems and develop an adaptive training procedure. The modified KRnet is referred to as the temporal KRnet (tKRnet) in this paper. The main contributions of this work are as follows. First, the basic layers in KRnet (where the temporal variable is not included) are systematically extended to be time-dependent, and the initial condition of the underlying stochastic dynamical system is encoded in the tKRnet as a prior distribution. Second, an adaptive training procedure for tKRnet is proposed. It is known that choosing proper collocation points is crucial for solving PDEs with deep learning based methods  \cite{TANG2022111080,tang2023daspinn}. To result in effective collocation points, our adaptive procedure has the following two main steps: training a tKRnet to approximate the solution of the Liouville equation, and using the trained tKRnet to generate collocation points for the next iteration. Through this procedure, the distribution of the collocation points becomes more consistent with the solution PDF after each iteration. Third, for the challenging problem of long-time integration associated with the Liouville equation, a temporal decomposition method is proposed, which provides guidance for applying tKRnet in this challenging problem. Lastly, a theoretical analysis is conducted to build the control of Kullback-Leibler (KL) divergence between the exact solution and the tKRnet approximation. 

The rest of the paper is organized as follows. Preliminaries of stochastic dynamical systems and the corresponding  Liouville equations are presented in \cref{sec:setup}. Detailed formulations of tKRnet are given in  \cref{sec:temporal-KRnet}. In \cref{sec:method}, our adaptive training procedure of tKRnet for solving the Liouville equation is presented, the corresponding temporal decomposition is discussed, and the analysis for the KL divergence between the exact solution and the tKRnet approximation is conducted.  Numerical results are discussed in \cref{sec:experiments}.  \Cref{sec:conclusions} concludes the paper. 

\section{Problem setup and preliminaries}
\label{sec:setup}
Let $\omega$ denote a random event, and $\bxi(\omega)\in\mathbb{R}^m$ and $\by_0(\omega)\in\mathbb{R}^n$ refer to random vectors, where $m$ and $n$ are positive integers. 
We consider the following stochastic dynamical system 
\begin{equation}
    \label{eq:stochastic-system}
    \frac{\dif \by(t)}{\dif t}= \mathg(\by,\bxi(\omega),t),\quad \by(0) = \by_0(\omega),
\end{equation}
where $t\in \mathI:=(0,T]$ ($T>0$ is a given final time),  $\by(t):=[y_1(t),\ldots,y_n(t)]^\top \in \mathbb{R}^n$ is a multi-dimensional stochastic process and  $\mathg:\mathbb{R}^n\times\mathbb{R}^{m}\times \mathI\rightarrow \mathbb{R}^n$ is a locally Lipschitz continuous 
vector function in terms of $\by$.  Let $\bx(t):=[\by(t),\bxi]^\top\in\mathbb{R}^{d}$ with $d=n+m$. The system \eqref{eq:stochastic-system} can be reformulated as 
\begin{equation}
    \label{eq:uncertain-para-init}
    \frac{\dif \bx(t)}{\dif t} = \mathf(\bx,t),\quad \bx(0) := [\by_0(\omega),\bxi(\omega)]^\top \textrm{ and } \mathf(\bx,t):=
    \left[
    \begin{aligned}
    \mathg(&\bx,t)
    \\
    &\mathbf{0}
    \end{aligned}
    \right].
\end{equation}
Our objective is to approximate the time-varying PDF of the state $\bx(t)$  
and the associated statistics.  

The PDF of $\bx(t)$, denoted as $p(\bx,t):\mathbb{R}^d\times \mathI\rightarrow \mathbb{R}_+$, satisfies the following Liouville equation 
\begin{equation}
    \label{eq:Liouville}
    \frac{\partial p(\bx,t)}{\partial t} + \nabla_{\bx} \cdot(p(\bx,t) \mathf(\bx,t)) = 0,\quad p(\bx,0) =p_0(\bx):= p_{\by}(\by,0)p_{\bxi}(\bxi), 
\end{equation}
where  $\nabla_{\bx} \cdot$ denotes the divergence in terms of $\bx$, $p_{\by}(\by,0)$ is the PDF of $\by_0(\omega)$, and $p_{\bxi}(\bxi)$ is the PDF of $\bxi(\omega)$.  
To improve the numerical stability, the logarithmic Liouville equation is proposed in 
\cite{villani2009Optimal,ben-hamu2022matching}, which can be written as \begin{equation}
    \label{eq:log-Liouville}
    \frac{\partial \log p(\bx,t)}{\partial t} +(\nabla_{\bx}\log p(\bx,t))\cdot \mathf(\bx,t)+\nabla_{\bx}\cdot \mathf(\bx,t)= 0.
\end{equation}
As $p(\bx,t)$ is a PDF for each $t\in \mathI$, it is required that \[\int_{\mathbb{R}^d} p(\bx,t) \dif \bx \equiv 1 \textrm{ and } p(\bx,t) \geq 0.\]
In addition, for any time $t$, the boundary condition for $p(\bx,t)$ is 
\[
p(\bx,t) \rightarrow 0 \textrm{ as } \|\bx\|_2\rightarrow \infty,\] 
where $\|\cdot\|_2$ indicates the $\ell_2$ norm.

\section{Temporal KRnet (tKRnet)}
\label{sec:temporal-KRnet}
The KRnet \cite{tang2020deep} is a flow-based generative model for density estimation or approximation, and its adaptive version is developed to solve the steady state Fokker-Planck equation in \cite{tang2023daspinn}. We systematically generalize the KRnet to a time-dependent setting in this section, which is referred to as the temporal KRnet (tKRnet) in the following, and develop an efficient adaptive training procedure for tKRnet to solve the Liouville equation in the next section.

Let $\bX\in\mathbb{R}^d$ be a random vector, which is associated with a time-dependent PDF $p_{\bX}(\bx,t)$. 
In this study, $p_{\bX}(\bx,t)$ is used to model the solution of  \eqref{eq:Liouville} (or \eqref{eq:log-Liouville}).
Let $\bZ\in\mathbb{R}^d$ be a random vector associated with a PDF $p_{\bZ}(\bz)$, where $p_{\bZ}(\bz)$ is a prior distribution (e.g., a Gaussian distribution). 
The main idea of time-dependent normalizing flows is to seek a time-dependent invertible mapping $\mathT:\mathbb{R}^d\times\mathbb{R}_+\rightarrow \mathbb{R}^d$ (i.e., $\bz=\mathT(\bx,t)$), and 
by the change of variables, the PDF $p_{\bX}(\bx,t)$ 
is given by 
\begin{equation}
\label{eq:time-change-variable}
     p_{\bX}(\bx,t) = p_{\bZ}(\bz)|\det \nabla_{\bx} \mathT(\bx,t)|, \textrm{ where } \bz = \mathT(\bx,t).
\end{equation}
Once the prior distribution $p_{\bZ}$ is specified, the explicit PDF for any random vector $\bX$ and time $t$ can be obtained through \eqref{eq:time-change-variable}. Additionally, exact random samples from $p_{\bX}(\bx,t)$ can be  obtained using the samples of $\bZ$ (from the prior) and the inverse of $\mathT$, i.e., $\bX=\mathT^{-1}(\bZ,t)$. In the rest of this paper, we let $\mathT$ denote our tKRnet, which is constructed by a sequence of time-dependent bijections. These include affine coupling layers, scale-bias, squeezing and nonlinear layers, which are defined as follows. 

\subsection{Time-dependent affine coupling layers}
\label{subsec:new-time-affine-coupling-layer}
A major part of tKRnet is affine coupling layers.  
Let $\bx=[\bx^{(1)},\bx^{(2)}]^\top$ be a partition of $\bx\in \mathbb{R}^d$ with $\bx^{(1)}\in\mathbb{R}^k$ and $\bx^{(2)}\in\mathbb{R}^{d-k}$,  where $k<d$ is a positive integer. For $t\in (0,T]$, 
the output of a time-dependent affine coupling layer $\tilde{\bx}=[\tilde{\bx}^{(1)},\tilde{\bx}^{(2)}]^\top$ is defined as  
\begin{equation}
\label{eq:affine-coupling}
    \begin{aligned}
    \tilde{\bx}^{(1)} &= \bx^{(1)}\\      \tilde{\bx}^{(2)}&= \bx^{(2)} + \frac{t}{T}\left(\alpha \bx^{(2)}\odot \tanh(\maths(\bx^{(1)}, t)) + e^{\bbeta}\odot\tanh(\matht(\bx^{(1)}, t))\right), 
    \end{aligned}
\end{equation}
where $0<\alpha<1$ is a fixed hyperparameter, $\bbeta\in\mathbb{R}^{d-k}$ is a trainable parameter, 
$\tanh$ is the hyperbolic tangent function, and $\odot$ is the Hadamard product or elementwise multiplication. We typically set $\alpha=0.6$ and  $k=\lfloor d/2 \rfloor$, where $\lfloor a \rfloor$ ($a \in \mathbb{R}$) is the largest integer that is smaller or equal to $a$. In \eqref{eq:affine-coupling}, $\maths:\mathbb{R}^{k+1}\rightarrow \mathbb{R}^{d-k}$ and $\matht:\mathbb{R}^{k+1}\rightarrow \mathbb{R}^{d-k}$ stand for scale and translation respectively,  
which are modeled by a neural network, i.e.,  \begin{align*}
    \left [\begin{aligned}
        \maths \\ \matht
    \end{aligned}\right] = \mathN(\bx^{(1)},t),
\end{align*}
where $\mathN:\mathbb{R}^{k+1}\rightarrow \mathbb{R}^{2(d-k)}$ is a feedforward neural network. 
Inspired by the work \cite{wang2021eigenvector}, we apply a random Fourier feature transformation before fully connected layers to alleviate possible challenges from multiscale problems. The neural network $\mathN$ is defined as 
\begin{equation}
\label{eq:affine-coupling-structure}
\begin{aligned}
    \bh_0 = &
    [\bx^{(1)},t]^\top
    ,\\
    \bh_1 = &\left[
    \begin{array}{c}
        \sin\left(\frac{1}{e^{\sigma}}\mathbf{F}\bh_0 + \bb_0\right)\\
    \cos\left(\frac{1}{e^{\sigma}}\mathbf{F}\bh_0 + \bb_0\right)\\
    \bh_0
    \end{array}
    \right],\\
    \bh_i = &\SiLU(\bW_{i-1}\bh_{i-1}+\bb_{i-1}),\quad \text{for}\quad  i=2,3,\ldots,M,\\
    \left[
    \begin{aligned}
    \maths\\
    \matht
    \end{aligned}
    \right]=&\bW_{M}\bh_M + \bb_{M},
\end{aligned}
\end{equation} 
where $\bh_0 \in \mathbb{R}^{k+1}$, $\bh_1\in \mathbb{R}^{d_h+k+1}$, and $\sigma\in \mathbb{R}$ is a learnable parameter. In \eqref{eq:affine-coupling-structure}, each entry of $\mathbf{F}\in \mathbb{R}^{(d_h/2)\times (k+1)}$ is sampled from the Gaussian distribution $\mathcal{N}(0, 1)$, and each entry of $\bb_0\in \mathbb{R}^{d_h/2}$ is sampled from the uniform distribution $\mathcal{U}(0,2\pi)$; $\bW_1\in\mathbb{R}^{d_h\times (d_h+k+1)}$ and  $\bb_1\in\mathbb{R}^{d_h}$ are the weight and the bias in the first fully connected layer, and $\bW_i\in\mathbb{R}^{d_h\times d_h}$ and $\bb_i\in \mathbb{R}^{d_h}$ are weight and bias coefficients in the $i$-th fully connected layer, for $i=2,\ldots,M-1$; $\bh_i\in\mathbb{R}^{d_h}$ is the hidden feature output by the $(i-1)$-th fully connected layer, and $\bW_M\in \mathbb{R}^{2(d-k)\times d_h}$ and $\bb_M\in\mathbb{R}^{2(d-k)}$ are weight and bias parameters in the output layer. Here, the trainable parameters in an affine coupling layer are $\{\bW_i, \bb_i| i=1,\ldots,M\}$. 
The sigmoid linear unit (SiLU) function \cite{hendrycks2016gaussian} is applied as the nonlinear activation function in the neural network, which is defined as \[\SiLU(x)=\frac{x}{1-e^{-x}}.\]

The Jacobian of $\tilde{\bx}$ at time $t$ with respect to $\bx$ can be obtained by 
\begin{equation}
    \label{eq:affine-coupling-jac}
    \nabla_{\bx}\tilde{\bx} = 
        \left[ \begin{aligned}\frac{\partial \tilde{\bx}^{(1)}}{\partial \bx^{(1)}} &\quad \frac{\partial \tilde{\bx}^{(1)}}{\partial \bx^{(2)}}\\
        \frac{\partial \tilde{\bx}^{(2)}}{\partial \bx^{(1)}}&\quad \frac{\partial \tilde{\bx}^{(2)}}{\partial \bx^{(2)}}\end{aligned}\right]=\left[\begin{aligned}
        \mathbb{I}_{k}&\quad \mathbf{0}\\
        \frac{\partial \tilde{\bx}^{(2)}}{\partial \bx^{(1)}} &\quad \diag(\mathbf{1}_{d-k}+\alpha\frac{t}{T}\tanh(\maths(\bx^{(1)}, t)))
        \end{aligned}\right],
\end{equation}
whose determinant can be efficiently computed due to the lower-triangular structure. 

It should be noted that, in the affine coupling layer, only $\bx^{(2)}$ is transformed to $\tilde{\bx}^{(2)}$, while $\bx^{(1)}$ remains unchanged. 
To completely update all components of $\bx$, the unchanged part $\bx^{(1)}$ is updated in the next affine coupling layer, 
\begin{align*}
    \tilde{\bx}^{(1)}&= \bx^{(1)} + \frac{t}{T}\left(\alpha \bx^{(1)}\odot \tanh(\maths(\bx^{(2)}, t)) + e^{\bbeta}\odot\tanh(\matht(\bx^{(2)}, t))\right)\\
    \tilde{\bx}^{(2)} &= \bx^{(2)}.
    \end{align*}

\subsection{Scale-bias, squeezing and nonlinear layers}
\label{subsec:layers-in-tKRnet}
We generalize the \textit{scale-bias layer} introduced in \cite{kingma2018glow} to a time-dependent setting. For an input $\bx\in\mathbb{R}^d$, the output of the time-dependent scale-bias layer is defined as  
\begin{equation}
    \label{eq:scale-bias}
    \tilde{\bx} = e^{\tanh(\bphi t)\odot \ba}\odot\bx + \tanh(\bphi t)\odot\bb, 
\end{equation}
where $\ba\in \mathbb{R}^d$, $\bb\in\mathbb{R}^d$, and $\bphi=[\phi_1,\ldots,\phi_d]^\top\in\mathbb{R}^d$ are trainable parameters, and each $\phi_i$ (for $i=1,\ldots,d$) is required to be positive.

\textit{Squeezing layer} $L_{S}$ is  applied to deactivate some dimensions through a mask   
\begin{equation}
    \label{eq:squeeze}
    \mathbf{m}=[ \underbrace{1,\ldots,1,}_{k}\underbrace{0,\ldots,0}_{d-k}]^{\top},
\end{equation}
where the components $\mathbf{m}\odot \bx$ are updated and the rest components $(1-\mathbf{m})\odot\bx$ are fixed from then on. 

\textit{Time-dependent nonlinear layer} $L_{N}$ extends the nonlinear layers introduced in \cite{tang2020deep} such that it is consistent with an identity transformation at time $0$. 
Let the interval $[0,1]$ be discretized by a mesh $0=s_0<s_1<\cdots<s_{\hat{m}+1}=1$ with the element size $h_i=s_{i+1}-s_i$, where $\hat{m}$ is a given positive integer. A nonlinear mapping $F(\cdot):[0,1]\rightarrow [0,1]$ is given by a quadratic function 
\begin{eqnarray}
\quad
F(s)=\frac{w_{i+1}-w_i}{2 h_i}(s-s_i)^2+w_i(s-s_i)+\sum_{k=0}^{i-1}\frac{w_k+w_{k+1}}{2}h_k,\, s\in[s_i,s_{i+1}].\label{eq:nonlinearF}
\end{eqnarray}
In \eqref{eq:nonlinearF}, the parameters $w_i$ are defined as $w_i={\hat{w}_i}/{c_w}$
with
\[
\hat{w}_i = e^{\tanh(\varphi_i t)\psi_i},\quad 
c_w = \sum_{i=0}^{\hat{m}} \frac{\hat{w}_i+\hat{w}_{i+1}}{2}h_i
,\] 
where $\varphi_i\in\mathbb{R}_+$ and $\psi_i\in\mathbb{R}$ are trainable parameters for $i=0,\ldots,\hat{m}$.  
As the support of $\bx$ in each dimension is $(-\infty,\infty)$,  
a mapping $\hat{F}(\cdot):(-\infty,\infty)\rightarrow(-\infty,\infty)$ is defined as
\begin{equation}
    \label{eq:non-linear-mapping}
    \hat{F}(x)=\begin{cases}\beta_s(x+a)-a,\quad x\in(-\infty,-a)\\
2aF(\frac{x+a}{2a})-a,\quad x\in[-a,a]\\
\beta_s(x-a)+a,\quad x\in(a,\infty),\\
\end{cases}
\end{equation}
where $\beta_s\in\mathbb{R}_+$ is a scaling factor and $a\in\mathbb{R}_+$ is a fixed hyperparameter. 
The inverse and the derivative with respect to $x$ of the nonlinear mapping can be explicitly computed \cite{TANG2022111080}.

\subsection{The overall structure of tKRnet}
\label{subsec:overall-structure-tKRnet}
\begin{figure}[htp!]
    \centering
    \includegraphics[scale=0.7]{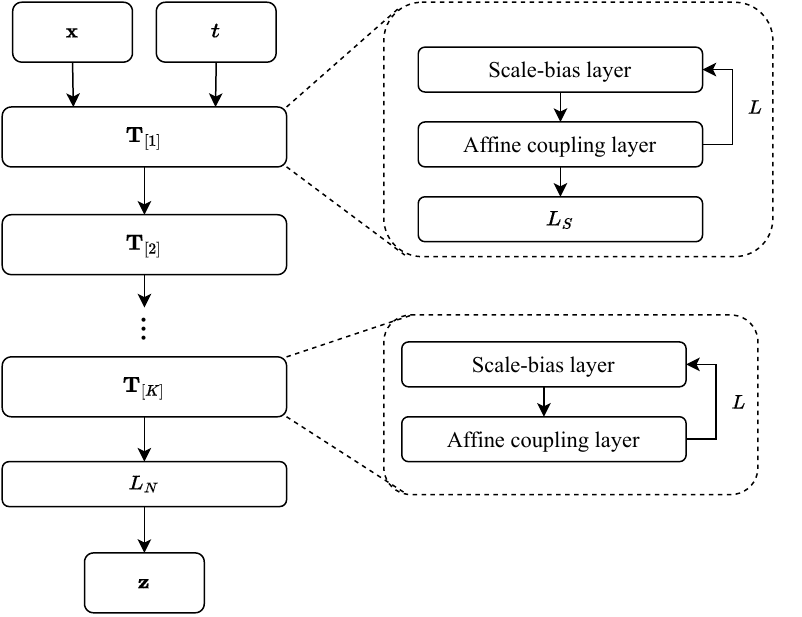}
    \caption{Overall structure of tKRnet.}
    \label{fig:overall-structure-temporal-KRnet}
\end{figure} 

Let $\bx=[\bx^{(1)},\ldots,\bx^{(K)}]^\top$ be a partition of $\bx$, where $\bx^{(i)}=[\bx^{(i)}_1,\ldots,\bx^{(i)}_m]^\top$ with $1\leq K\leq d$, $1\leq m\leq d$, and $\sum_{i=1}^K \operatorname{dim}(\bx^{(i)}) = d$.  
The partition for $\bz$ is the same as that for $\bx$,  
and the inverse transformation of the tKRnet is set to a lower-triangular transport map, i.e., 
\begin{align*}
    \bx = \mathT^{-1}(\bz,t)=\left[\begin{aligned}
    &\mathT_{K}^{-1}(\bz^{(1)},t)\\
    &\mathT_{K-1}^{-1}(\bz^{(1)},\bz^{(2)},t)\\
     &\vdots\\
    &\mathT_{1}^{-1}(\bz^{(1)},\ldots,\bz^{(K)},t)\end{aligned} \right],
\end{align*}
where $\mathT_i$ are time-dependent invertible mappings for $i=1,\ldots,K$. The forward computation of the tKRnet is defined as a composite form  
\begin{eqnarray}
\bz=\mathT(\bx,t;\Theta)=L_N\circ\mathT_{[K]}\circ\cdots\circ\mathT_{[1]}(\bx,t), \label{eq:Ttheta}
\end{eqnarray}
where $\Theta$ refer to the parameters in this transformation, and $\mathT_{[i]}$ is defined as
\begin{eqnarray}
\mathT_{[i]}:=\begin{cases}
L_S\circ\mathT_{[i,L]}\circ\cdots\circ\mathT_{[i,1]},\quad &i=1,\ldots,K-1,\\
\mathT_{[i,L]}\circ\cdots\circ\mathT_{[i,1]},\quad &i=K.
\end{cases}\label{eq:Ti}
\end{eqnarray}
In \eqref{eq:Ttheta}--\eqref{eq:Ti}, the invertible mapping $\mathT_{[\cdot,\cdot]}$ is composed by one scale-bias layer \eqref{eq:scale-bias} and one affine coupling layer (see section \ref{subsec:new-time-affine-coupling-layer}), and $L_N,L_S$ indicate the nonlinear layer \eqref{eq:nonlinearF} and the squeezing layer \eqref{eq:squeeze} respectively. 
The overall structure of tKRnet is illustrated in \cref{fig:overall-structure-temporal-KRnet}.

Without loss of generality, denoting $\bx_{[i]}=\mathT_{[i]}(\bx_{[i-1]},t)$, the partitions $\bx_{[i]}=[\bx_{[i]}^{(1)},\ldots,\bx_{[i]}^{(K)}]^\top$ for $i=1,\ldots,K$ are set to the same.  
At the beginning, $\mathT_{[1]}$ is applied to the partition $\bx_{[0]}=[\bx_{[0]}^{(1:K-1)},\bx_{[0]}^{(K)}]^\top$, where $\bx_{[0]}^{(1:K-1)}$ includes $\{\bx_{[0]}^{(1)},\ldots,\bx_{[0]}^{(K-1)}\}$. After that, ${\bx}_{[1]}=\mathT_{[1]}({\bx}_{[0]},t)$, 
and from then on, the last partition ${\bx}_{[1]}^{(K)}$ remains fixed, i.e., ${\bx}_{[i]}^{(K)}={\bx}_{[1]}^{(K)}$ for $i>1$.  
In general, after applying the transformation $\mathT_{[i]}$, the $(K-i+1)$-th partition of ${\bx}_{[i]}$, i.e., ${\bx}_{[i]}^{(K-i+1)}$, is  deactivated in addition to the dimensions fixed at previous stages.    

Based on the tKRnet $\mathT(\cdot,\cdot;\Theta)$, the PDF $p_{\Theta}(\bx,t)$ can be obtained by \eqref{eq:time-change-variable} and the chain rule, 
\begin{equation}
    \begin{aligned}
        p_{\Theta}(\bx,t) &= p_{\bZ}(\mathT(\bx,t;\Theta))|\det \nabla_{\bx} \mathT(\bx,t;\Theta)|\\
        &=p_{\bZ}(\mathT(\bx,t;\Theta))|\det\nabla_{{\bx}_{[K]}}L_N ({\bx}_{[K]},t)|\prod_{i=1}^{K}|\det\nabla_{{\bx}_{[i-1]}}\mathT_{[i]}({\bx}_{[i-1]},t)|,
    \end{aligned}
    \label{eq:pdf-change}
\end{equation}
where $\bx_{[0]}=\bx$, $\bx_{[i]} = \mathT_{[i]}(\bx_{[i-1]},t)$, for $i=1,\ldots,K$.
To use $p_{\Theta}(\bx,t)$ as a PDF model for the solution of the Liouville equation \eqref{eq:Liouville}, the prior distribution $p_{\bZ}$ is set to the initial PDF $p_0$. At time $t=0$, all sublayers of tKRnet are identity mappings, and $p_{\Theta}(\bx,0)=p_0(\mathT(\bx,0;\Theta))|\det\nabla_{\bx}\mathT(\bx,0;\Theta)|=p_0(\bx)$.

\section{Adaptive sampling based physics-informed training for density approximation}
\label{sec:method}
In this section, tKRnet is applied to solve the time-dependent Liouville equation \eqref{eq:Liouville} or its variant  \eqref{eq:log-Liouville}. We pay particular attention to adaptive sampling and long-time integration.

\subsection{Physics-informed training and adaptive sampling}
\label{subsec:adaptive-training}
Adaptivity is crucial for numerically solving PDF equations. The work \cite{Cho2013adaptive} uses spectral elements on an adaptive non-conforming grid to discretize the spatial domain of the PDF equation and track the time-dependent PDF support. The work \cite{TANG2022111080} shows that adaptive sampling strategy can help normalizing flow models effectively approximate solutions to steady Fokker-Planck equations. We propose a physics-informed method consisting of multiple adaptivity iterations. Each adaptivity iteration has two steps: 1) training tKRnets by minimizing the total PDE residual on collocation points in the training set; 2) generating collocation points dynamically through the trained tKRnet to update the training set. 

Let $p_{\Theta}(\bx,t)$ be the approximate PDF induced by tKRnet. 
The residual given by the Liouville equation \eqref{eq:Liouville} 
is defined as 
\begin{equation}
    \label{eq:residual-error}
    r(\bx,t;\Theta) := \frac{\partial p_{\Theta}(\bx,t)}{\partial t}+\nabla_{\bx}\cdot (p_{\Theta}(\bx,t)\mathf(\bx,t)).
\end{equation} 
In high-dimensional problems, the value of $p_{\Theta}(\bx,t)$ may be too small, which results in underflow such that the residual cannot be effectively minimized.
To alleviate this issue, the logarithmic Liouville equation \eqref{eq:log-Liouville} is considered, and the corresponding residual is defined as 
\begin{equation}
    \label{eq:log-residual-error}
r_{\log}(\bx,t;\Theta):=\frac{\partial \log p_{\Theta}(\bx,t)}{\partial t} +(\nabla_{\bx}\log p_{\Theta}(\bx,t))\cdot\mathf(\bx,t) +\nabla_{\bx}\cdot \mathf(\bx,t).
\end{equation}
It is clear that $r(\bx,t;\Theta)=p_{\Theta}(\bx,t)r_{\log}(\bx,t;\Theta)$. 
Letting $p_{c}(\bx,t)$ be a PDF in the space-time domain, the following loss functional is defined, 
\begin{equation}
    \mathcal{L}(p_{\Theta}(\bx,t))
    :=\mathbb{E}_{p_{c}(\bx,t)}
    (|r_{\log}(\bx,t;\Theta)|^2),
\end{equation}
where $\mathbb{E}_{p_{c}}$ refers to the expectation with respect to $p_{c}(\bx,t)$.  
A set of collocation points $\mathcal{S}=\{(\resx^{(i)},\rest^{(i)})\}_{i=1}^{N_r}$ are drawn from $p_{c}(\bx,t)$, and these points are employed to estimate the loss function as 
\begin{equation}
    \label{eq:loss-func}
   \mathcal{L}(p_{\Theta}(\bx,t)) \approx\hat{\mathcal{L}}(p_{\Theta}(\bx,t)):=\frac{1}{N_r}\sum_{i=1}^{N_r} |r_{\log}(\resx^{(i)},\rest^{(i)};\Theta)|^2.
\end{equation}
Optimal parameters for $p_{\Theta}(\bx,t)$ are chosen by minimizing the approximate loss function, i.e.,
\begin{eqnarray}
\Theta^*=\argmin_{\Theta} \hat{\mathcal{L}}(p_{\Theta}(\bx,t)).
\label{eq:opt}
\end{eqnarray}
A variant of the stochastic gradient descent method, AdamW \cite{loshchilov2018decoupled}, is applied to solve the optimization problem \eqref{eq:opt}. The learning rate scheduler for the optimizer is the cosine annealing method \cite{loshchilov2017sgdr}. 
More specifically, the set $\mathcal{S}$ is divided into $N_b$ mini-batches $\{\mathcal{S}^{(n)}\}_{n=1}^{N_b}$. Let $\Theta_{e,n}$ be the model parameters at the $n$-th step of epoch $e$ with $e=1,\ldots,N_E$ and $n=1,\ldots,N_b$.  
The model parameters are updated as \begin{equation}
    \label{eq:SGD-update-parameters}
    \Theta_{e,n} = \Theta_{e,n-1} - \tau \nabla_{\Theta_{e,n-1}} \frac{1}{|\mathcal{S}^{(n)}|} \sum_{(\resx,\rest)\in \mathcal{S}^{(n)}} |r_{\log}(\resx,\rest;\Theta_{e,n-1})|^2,
\end{equation}
where $\tau>0$ is the learning rate. After each optimization step, the learning rate is decreased by the learning rate scheduler.

To enhance the accuracy of the final approximation, we propose the following adaptive sampling strategy. 
Let $N_r=MJ$, where $M$ and $J$ are positive integers. 
On $(0,T]$, $\rest^{(i)}$ (for $i=1,\ldots,N_r$) are set to 
\begin{equation}
    \label{eq:residual-time-points}
    \rest^{(i)}= \lceil i/M \rceil \Delta t,
\end{equation}
where $\Delta t={T}/{(J-1)}$ and $\lceil i/M \rceil$ is the smallest integer that is larger or equal to $i/M$. Letting $\mathcal{S}_{\txttime}=\{\rest^{(i)}\}_{i=1}^{N_r}$, the training set is initialized as $\mathcal{S}_0:=\{(\resx^{(0,i)},\rest^{(i)})\}_{i=1}^{N_r}$, where
$\{\resx^{(0,i)}\}_{i=1}^{N_r}$ are the initial spatial collocation points (e.g., samples of a uniform distribution).  
The set $\mathcal{S}_0$ is divided into $N_b$ mini-batches $\{\mathcal{S}_0^{(n)}\}_{n=1}^{N_b}$. 
The tKRnet (introduced in  section \ref{subsec:overall-structure-tKRnet}) is initialized as $\mathT(\bx,t;\Theta^{(0)})$. 
In general, the $i$-th ($i=1,\ldots,N_r$) spatial collocation point generated at the $k$-th adaptivity iteration step is denoted by $\resx^{(k,i)}$. 
The parameters at the $k$-th adaptivity iteration step, the $e$-th epoch and the $n$-th optimization iteration step is denoted by $\Theta_{e,n}^{(k)}$. 

Starting with $\Theta_{1,0}^{(1)}=\Theta^{(0)}$, the tKRnet is trained through solving \eqref{eq:opt} with collocation points $\mathcal{S}_0$,  
and the trained tKRnet at adaptivity iteration step one is denoted by  $\mathT(\bx,t;\Theta^{*,(1)})$, where $\Theta^{*,(1)}$ are the parameters of the trained tKRnet at this step. The PDF  $p_{\Theta^{*,(1)}}=p_0(\mathT(\bx,t;\Theta^{*,(1)}))|\det \nabla_{\bx} \mathT(\bx,t;\Theta^{*,(1)})|$ (see \eqref{eq:pdf-change}) is then used to generate new collocation points  
$\mathcal{S}_1 =\{(\resx^{(1,i)},\rest^{(i)})\}_{i=1}^{N_r}$, where $\rest^{(i)}\in \mathcal{S}_{\txttime}$ remains unchanged and each $\resx^{(1,i)}$ is drawn from $p_{\Theta^{*,(1)}}(\bx,\rest^{(i)})$. To generate each $\resx^{(1,i)}$, a sample point $\bz^{(1,i)}$ is first generated using $p_{\bZ}=p_0$, and then  $\resx^{(1,i)}=\mathT^{-1}(\bz^{(1,i)},\rest^{(i)};\Theta^{*,(1)})$. 
Next, starting with $\mathT(\bx,t;\Theta^{(2)}_{1,0})$ where $\Theta^{(2)}_{1,0}:=\Theta^{*,(1)}$, we continue the training process with the collocation points $\mathcal{S}_1$ to obtain $\mathT(\bx,t;\Theta^{*,(2)})$.  
In general, at adaptivity iteration step $k$, the collocation points  $\mathcal{S}_{k}$ are generated using $p_{\Theta^{*,(k)}}(\bx,t)$, and the tKRnet with the initial parameters $\Theta^{(k+1)}_{1,0}:=\Theta^{*,(k)}$ is subsequently updated to $\mathT(\bx,t;\Theta^{*,(k+1)})$ through solving \eqref{eq:opt}. 
The adaptivity iteration continues until the maximum number of iterations is achieved, which is denoted by $N_{\adaptive}\in\mathbb{N}_+$. This adaptive procedure is summarized in \cref{alg:adaptive-sample-flow}. 
\begin{algorithm}[!htp]
\caption{Adaptive sampling based physics-informed training for tKRnet}
\label{alg:adaptive-sample-flow}
\algorithmicrequire{ The initial tKRnet $\mathT(\bx,t;\Theta^{(0)})$, the number of collocation points $N_r$, the number of epochs $N_{E}$, the maximum number of adaptivity iterations $N_{\adaptive}$, the learning rate $\tau$, and the number of mini-batches $N_b$.}
\begin{algorithmic}[1]
\STATE{Generate time points  $\mathcal{S}_{\txttime}=\{\rest^{(i)}\}_{i=1}^{N_r}$ (see \eqref{eq:residual-time-points}).}
\STATE{Draw initial samples  $\{\resx^{(0,i)}\}_{i=1}^{N_r}$ from a given distribution, 
and let $\mathcal{S}_0=\{(\resx^{(0,i)},\rest^{(i)})\}_{i=1}^{N_r}$.}
\STATE{Divide $\mathcal{S}_0$ into $N_b$ mini-batches $\{\mathcal{S}^{(n)}_{0}\}_{n=1}^{N_b}$.}
\FOR{$k=1$,\ldots,$N_{\adaptive}$}
\IF{$k=1$}
\STATE{Let $\Theta^{(k)}_{1,0}=\Theta^{(0)}$.}
\ELSE
\STATE{Let $\Theta^{(k)}_{1,0}=\Theta^{*,(k-1)}$.}
\ENDIF
\STATE{Initialize the AdamW optimizer and the cosine annealing learning rate scheduler.}
\FOR{$e=1$,\ldots,$N_E$}
\FOR{$n=1$,\ldots,$N_b$}
\STATE{Compute the loss $\hat{\mathcal{L}}(p_{\Theta^{(k)}_{e,n-1}}(\bx,t))$ (see \eqref{eq:loss-func}) on the mini-batch $\mathcal{S}^{(n)}_{k-1}$.}
\STATE{Update $\Theta^{(k)}_{e,n}$ using the AdamW optimizer with the learning rate $\tau$ (see \eqref{eq:SGD-update-parameters}).}
\STATE{The learning rate scheduler decreases the learning rate.}
\ENDFOR
\IF{$e<N_E$}
\STATE{Let $\Theta^{(k)}_{e+1,0}=\Theta^{(k)}_{e,N_b}$.}
\STATE{Shuffle the mini-batches $\{\mathcal{S}^{(n)}_{k-1}\}_{n=1}^{N_b}$ of the set $\mathcal{S}_{k-1}$.}
\ENDIF
\ENDFOR
\STATE{Let $\Theta^{*,(k)}=\Theta^{(k)}_{N_E,N_b}$.}
\STATE{Draw $\resx^{(k,i)}$ from $p_{\Theta^{*,(k)}}(\bx,\rest^{(i)})$, and let $\mathcal{S}_k=\{(\resx^{(k,i)},\rest^{(i)})\}_{i=1}^{N_r}$.}
\STATE{Divide $\mathcal{S}_k$ into $N_b$ mini-batches $\{\mathcal{S}_{k}^{(n)}\}_{n=1}^{N_b}$.}

\ENDFOR
\STATE{Let $\Theta^*=\Theta^{*,(N_{\adaptive})}$.}
\end{algorithmic}
\algorithmicensure{ The trained tKRnet $\mathT(\bx,t;\Theta^*)$, and the approximate solution $p_{\Theta^*}(\bx,t)$.}
\end{algorithm}

\subsection{Temporal decomposition for long-time integration}
\label{subsec:temporal-dcomposition-long-time}
The performance of PINN may deteriorate for evolution equations when the time domain becomes large \cite{wang2023long-time}. \Cref{alg:adaptive-sample-flow} suffers a similar issue since it is defined in the framework of PINN.
Some remedies \cite{meng2020PPINN,penwarden2023temporaldecomposition} have been proposed to alleviate this issue. In this work, we employ a temporal decomposition method when implementing \cref{alg:adaptive-sample-flow} on a large time domain. We will demonstrate numerically that coupling temporal decomposition and adaptive sampling yields efficient performance for long-time integration, although other techniques \cite{wang2023long-time,penwarden2023temporaldecomposition} can also be employed for further refinement.

We decompose the time interval $(0,T]$ into $N_{sub}$ sub-intervals  $(T_{i-1},T_i]$ (for $i=1,\ldots,N_{sub}$), where $0=T_0<T_1<\cdots <T_{N_{sub}}=T$. 
For the $i$-th sub-interval $(T_{i-1},T_{i}]$, assuming that the PDF $p_{\Theta}(\cdot,T_{i-1})$ is given (e.g., for the interval $(T_{0},T_1]$, $p_{\Theta}(\cdot,T_{0})$ is given a priori), our goal is to train tKRnet $\mathT(\cdot,\cdot;\Theta^{(i,k)})$, where $\Theta^{(i,k)}$ includes the model parameters and $k$ denotes the adaptivity iteration step (see line 4 of \cref{alg:adaptive-sample-flow}). The trained parameters of the tKRnet for $i$-th sub-interval are denoted by $\Theta^{(i),*}$.  
The following two choices for temporal decomposition are considered in this work to train the local tKRnets. The first choice is to 
keep the same tKRnet structure for different sub-intervals (while the local tKRnet parameters are different), and to introduce the cross entropy to maintain the consistency between two adjacent sub-intervals; the second choice is to use $p_{\Theta}(\cdot,T_{i-1})$ as the prior distribution to construct a new tKRnet for each sub-interval $(T_{i-1},T_i]$. 

In the first choice, the solution of \eqref{eq:Liouville} or \eqref{eq:log-Liouville} on $(T_{i-1}, T_i]$ is approximated by  $p(\bx,t)\approx p_{\Theta^{(i,k)}}(\bx,t):=p_0(\bz)|\det\nabla_{\bx}\mathT(\bx,t;\Theta^{(i,k)})|$ where $\bz=\mathT(\bx,t;\Theta^{(i,k)})$ and $\mathT$ is introduced in \cref{sec:temporal-KRnet}. Here, we let $p_{\Theta^{(i,k)}}(\bx,T_{i-1})\approx p_{\Theta^{(i-1),*}}(\bx,T_{i-1})$, where $p_{\Theta^{(i-1),*}}(\bx,T_{i-1})$ is the given tKRnet solution at $t=T_{i-1}$. 
The cross entropy between $p_{\Theta^{(i-1),*}}(\bx,T_{i-1})$ and $p_{\Theta^{(i,k)}}(\bx,T_{i-1})$ is next introduced to maintain the consistency at time $t=T_{i-1}$. 
That is, the loss function on $(T_{i-1},T_i]$ is defined as 
\begin{equation}
\label{eq:temporal-decomposition-loss}
\begin{aligned}
\hat{\mathcal{L}}(p_{\Theta^{(i,k)}}(\bx,t))&=\frac{1}{N_{r}}\sum_{j=1}^{N_r}|r_{\log}(\resx^{(j)},\rest^{(j)};\Theta^{(i,k)})|^2 \\
&\quad - \frac{1}{N_{r}}\sum_{j=1}^{N_{r}}\log p_{\Theta^{(i,k)}}(\bx_{\txtinterface}^{(j)},T_{i-1}),
\end{aligned}
\end{equation}
where $\bx_{\txtinterface}^{(j)}$ are drawn from $p_{\Theta^{(i-1),*}}(\cdot,T_{i-1})$. Note that the loss function for $\mathT(\cdot,\cdot;\Theta^{(1,k)})$ (on the sub-interval $(T_0,T_1]$) is  
\begin{eqnarray}
\hat{\mathcal{L}}(p_{\Theta^{(1,k)}}(\bx,t))=\frac{1}{N_{r}}\sum_{j=1}^{N_r}|r_{\log}(\resx^{(j)},\rest^{(j)};\Theta^{(1,k)})|^2,
\label{eq:loss_sub1}
\end{eqnarray}
where only the residual is considered, because the initial distribution at $t=0$ is used as the prior distribution for the tKRnet. Then, \cref{alg:adaptive-sample-flow} is applied with \eqref{eq:loss_sub1} is obtain the local solution $p_{\Theta^{(1),*}}$ for the first sub-interval $(T_0,T_1]$. After that,  \cref{alg:adaptive-sample-flow} is applied with \eqref{eq:temporal-decomposition-loss} to train $p_{\Theta^{(2),*}}$, and this procedure is repeated for computing $(T_{i-1},T_i]$ for $i=3,\ldots,N_{sub}$.

In the second choice, the local tKRnet for each sub-interval is rebuilt with the prior distribution obtained at the previous sub-interval. 
The local tKRnet for $(T_{i-1}, T_i]$ is denoted by $\subT{i}(\bx, t; \Theta^{(i,k)})$.  
The tKRnet solution of \eqref{eq:Liouville} (or \eqref{eq:log-Liouville}) at $t=T_{i-1}$ in this choice is defined as, 
\begin{equation}
\label{eq:prior-density-computation}
\begin{aligned}
    p_{\Theta}(\bx,T_{i-1})&=p_0(\hat{\bx}_{0}) \prod_{j=1}^{i-1}|\det\nabla_{\hat{\bx}_{i-j}}\subT{i-j}(\hat{\bx}_{i-j},T_{i-j};\Theta^{(i-j),*})|, 
\end{aligned}    
\end{equation}
where 
\begin{equation}
\label{eq:prior-density-computation-1}
\begin{aligned}
    \hat{\bx}_{i-1}=\bx, \quad
\hat{\bx}_{i-j-1} &= \subT{i-j}(\hat{\bx}_{i-j},T_{i-j};\Theta^{(i-j),*}),\quad j=1,\ldots,i-1.
\end{aligned}
\end{equation}
In \eqref{eq:prior-density-computation}--\eqref{eq:prior-density-computation-1}, $\{\Theta^{(i-j),*}|j=1,\ldots,i-1\}$ denotes the set of trained parameters associated with local tKRnets for the previous sub-intervals. With $p_{\Theta}(\cdot, T_{i-1})$ defined in \eqref{eq:prior-density-computation}, the tKRnet solution for $t\in(T_{i-1},T_i]$ is obtained through 
\begin{equation}
\label{eq:prenets-as-prior-approximate-PDF}
\begin{aligned}
    p_{\Theta}(\bx,t)&=p_{\Theta}(\subT{i}(\bx,t;\Theta^{(i,k)}),T_{i-1})|\det \nabla_{\bx} \subT{i}(\bx,t;\Theta^{(i,k)})|,
\end{aligned}
\end{equation}
where $\Theta=\{\Theta^{(i,k)},\Theta^{(i-j),*}|j=1,\ldots,i-1\}$. Note that the trainable weights at this stage are $\Theta^{(i,k)}$ while the other parameters remain fixed, and $\subT{i}(\cdot, T_{i-1}; \Theta^{(i,k)})$ is an identity mapping, as  $p_{\Theta}(\cdot,T_{i-1})$ is the prior distribution. As the tKRnet $\subT{i}$ implemented in this method is dependent on $T_{i-1}$ and $t$, and the affine coupling layer (see \eqref{eq:affine-coupling}) is modified as 
\begin{equation}
    \label{eq:affine-coupling-temporal-decom}
    \begin{aligned}
    \tilde{\bx}^{(1)} &= \bx^{(1)}\\      \tilde{\bx}^{(2)}&= \bx^{(2)} + \frac{t-T_{i-1}}{T_{i}-T_{i-1}}\left(\alpha \bx^{(2)}\odot \tanh(\maths(\bx^{(1)}, t, T_{i-1})) + e^{\bbeta}\odot\tanh(\matht(\bx^{(1)}, t, T_{i-1}))\right), 
    \end{aligned}
\end{equation}
where $[\maths , \matht ]^\top = \mathN(\bx^{(1)},t-T_{i-1}, T_{i-1})$. In scale-bias layers and nonlinear layers (presented in \cref{sec:temporal-KRnet}), the temporal variable $t$ is replaced by  $t-T_{i-1}$ in this setting. For example, scale-bias layers here are defined as 
\begin{equation}
\label{eq:scale-bias-temporal-decom}
    \tilde{\bx} = e^{\tanh(\bphi (t-T_{i-1}))\odot \ba}\odot\bx + \tanh(\bphi (t-T_{i-1}))\odot\bb.
\end{equation}
After these slight modifications, \cref{alg:adaptive-sample-flow} can be implemented to obtain a local solution $p_{\Theta}(\bx,t)$ \eqref{eq:prenets-as-prior-approximate-PDF} on each sub-interval $(T_{i-1},T_i]$ for $i=1,\ldots,T_{sub}$.

\subsection{Theoretical properties}
\label{subsec:theory-analysis}
In this section, we show that the forward Kullback-Leibler (KL) divergence between the exact solution $p(\bx,t)$ of \eqref{eq:Liouville} (or \eqref{eq:log-Liouville}) and the tKRnet solution $p_{\Theta}(\bx,t)$ (see \eqref{eq:residual-error}) can be bounded by the residual $r_{\log}(\bx,t;\Theta)$ \eqref{eq:log-residual-error}. 

\begin{theorem}[Control of KL divergence]
\label{thm:time-derivative-KL}
Assume $p(\bx,t)\rightarrow 0, p_{\Theta}(\bx,t)\rightarrow 0$ as $\|\bx\|_2\rightarrow \infty$. Denote $\Omega_{C}=\{\bx|\|\bx\|_2\leq C\}$ with $C>0$. For any $t\in[0,T]$, assume $\lim_{C\rightarrow \infty} \oint_{\partial \Omega_{C}} \log\left(\frac{p(\bx,t)}{p_{\Theta}(\bx,t)}\right)p(\bx,t)\mathf(\bx,t) \cdot \boldsymbol{n} \dif \boldsymbol{s}=0$, where $\boldsymbol{n}$ is the outward pointing unit normal. 
The foward KL divergence between the exact solution $p(\bx,t)$ of \eqref{eq:Liouville} (or \eqref{eq:log-Liouville}) and the tKRnet solution $p_{\Theta}(\bx,t)$ (see \eqref{eq:residual-error}) satisfies 
\begin{align*}
    \frac{\dif }{\dif t} D_{KL}(p(\bx,t)||p_{\Theta}(\bx,t))&\leq \int_{\mathbb{R}^d} |r_{\log} (\bx,t;\Theta)|p(\bx,t)\dif \bx.
\end{align*}
\end{theorem}
\begin{proof}
For KL divergence, we have 
\begin{align*}
    \frac{\dif }{\dif t} &D_{KL}(p(\bx,t)||p_{\Theta}(\bx,t))\\
    &=\frac{\dif }{\dif t}\int_{\mathbb{R}^d} \log \left(\frac{p(\bx,t)}{p_{\Theta}(\bx,t)}\right)p(\bx,t) \dif \bx\\
    &=\int_{\mathbb{R}^d} \frac{\partial }{\partial t}\log \left(\frac{p(\bx,t)}{p_{\Theta}(\bx,t)}\right)p(\bx,t) \dif \bx+  \int_{\mathbb{R}^d}\log \left(\frac{p(\bx,t)}{p_{\Theta}(\bx,t)}\right)\frac{\partial p(\bx,t)}{\partial t}\dif \bx\\
    &=\int_{\mathbb{R}^d} \frac{\partial p(\bx,t)}{\partial t} \dif \bx - \int_{\mathbb{R}^d} \frac{p(\bx,t)}{p_{\Theta}(\bx,t)}\frac{\partial p_{\Theta}(\bx,t)}{\partial t} \dif \bx + \int_{\mathbb{R}^d} \log \left(\frac{p(\bx,t)}{p_{\Theta}(\bx,t)}\right)\frac{\partial p(\bx,t)}{\partial t}\dif \bx\\
    &= \frac{\partial \int_{\mathbb{R}^d} p(\bx,t)\dif \bx}{\partial t} - \int_{\mathbb{R}^d} \left(\frac{p(\bx,t)}{p_{\Theta}(\bx,t)}\right)(r(\bx,t;\Theta) - \nabla_{\bx}\cdot(p_{\Theta}(\bx,t)\mathf(\bx,t))) \dif \bx \\
    &\quad - \int_{\mathbb{R}^d} \log \left(\frac{p(\bx,t)}{p_{\Theta}(\bx,t)}\right)\nabla_{\bx}\cdot(p(\bx,t)\mathf(\bx,t))\dif \bx. 
\end{align*}
Since $\int_{\mathbb{R}^d}p(\bx,t) \dif \bx = 1$, 
\[\frac{\partial \int_{\mathbb{R}^d} p(\bx,t)\dif \bx}{\partial t}=\frac{\partial 1}{\partial t}=0.\]
Then, 
\begin{align*}
    \frac{\dif }{\dif t} D_{KL}(p(\bx,t)||p_{\Theta}(\bx,t))=&
    \underbrace{\int_{\mathbb{R}^d} \frac{p(\bx,t)}{p_{\Theta}(\bx,t)}\nabla_{\bx}\cdot(p_{\Theta}(\bx,t)\mathf(\bx,t)) \dif \bx}_{I_1}\\
    &- \underbrace{\int_{\mathbb{R}^d} \log\left (\frac{p(\bx,t)}{p_{\Theta}(\bx,t)}\right)\nabla_{\bx}\cdot(p(\bx,t)\mathf(\bx,t)) \dif \bx}_{I_2}\\
    &- \int_{\mathbb{R}^d} \frac{r(\bx,t;\Theta)}{p_{\Theta}(\bx,t)} p(\bx,t) \dif \bx.
\end{align*}
For $I_1$, integration by parts yields that
\begin{align*}
    \int_{\mathbb{R}^d} &\frac{p(\bx,t)}{p_{\Theta}(\bx,t)}\nabla_{\bx}\cdot(p_{\Theta}(\bx,t)\mathf(\bx,t)) \dif \bx\\
    &=\int_{\mathbb{R}^d} \nabla_{\bx}\cdot (p(\bx,t)\mathf(\bx,t)) \dif \bx- \int_{\mathbb{R}^d} \nabla_{\bx} \left(\frac{p(\bx,t)}{p_{\Theta}(\bx,t)}\right)\cdot \mathf(\bx,t) p_{\Theta}(\bx,t)\dif \bx\\
    &=- \int_{\mathbb{R}^d} \nabla_{\bx} \left(\frac{p(\bx,t)}{p_{\Theta}(\bx,t)}\right)\cdot \mathf(\bx,t) p_{\Theta}(\bx,t)\dif \bx,
\end{align*}
where the second equality is obtained using 
\begin{align*}
    \int_{\mathbb{R}^d} \nabla_{\bx}\cdot (p(\bx,t)\mathf(\bx,t)) \dif \bx =  -\int_{\mathbb{R}^d} \frac{\partial p(\bx,t)}{\partial t} \dif \bx=0.
\end{align*}
Similarly, $I_2$ can be rewritten as 
\begin{align*}
    \int_{\mathbb{R}^d}& \log\left(\frac{p(\bx,t)}{p_{\Theta}(\bx,t)}\right)\nabla_{\bx}\cdot(p(\bx,t)\mathf(\bx,t)) \dif \bx\\
    &=\int_{\mathbb{R}^d} \nabla_{\bx}\cdot \left(\log\left(\frac{p(\bx,t)}{p_{\Theta}(\bx,t)}\right)p(\bx,t)\mathf(\bx,t)\right)\dif \bx\\ 
    &\quad - \int_{\mathbb{R}^d} \nabla_{\bx}\log\left(\frac{p(\bx,t)}{p_{\Theta}(\bx,t)}\right)\cdot\mathf(\bx,t) p(\bx,t)\dif \bx\\
    &=\int_{\mathbb{R}^d} \nabla_{\bx}\cdot \left(\log\left(\frac{p(\bx,t)}{p_{\Theta}(\bx,t)}\right)p(\bx,t)\mathf(\bx,t)\right)\dif \bx\\
    &\quad- \int_{\mathbb{R}^d} \nabla_{\bx}\left(\frac{p(\bx,t)}{p_{\Theta}(\bx,t)}\right)\cdot\mathf(\bx,t) p_{\Theta}(\bx,t)\dif \bx.
\end{align*}
By the divergence theorem, 
\begin{align*}
    &\int_{\mathbb{R}^d} \nabla_{\bx}\cdot \left(\log\left(\frac{p(\bx,t)}{p_{\Theta}(\bx,t)}\right)p(\bx,t)\mathf(\bx,t)\right)\dif \bx\\
    &\quad=\lim_{C\rightarrow \infty}\oint_{\partial \Omega_{C}} \log\left(\frac{p(\bx,t)}{p_{\Theta}(\bx,t)}\right)p(\bx,t)\mathf(\bx,t) \cdot \boldsymbol{n} \dif \boldsymbol{s}=0.
\end{align*}
Finally, we get 
\begin{align*}
    \frac{\dif }{\dif t} D_{KL}(p(\bx,t)||p_{\Theta}(\bx,t))&=
     - \int_{\mathbb{R}^d} \frac{r(\bx,t;\Theta)}{p_{\Theta}(\bx,t)} p(\bx,t) \dif \bx\\
    &= - \int_{\mathbb{R}^d} r_{\log}(\bx,t;\Theta) p(\bx,t) \dif \bx\leq \int_{\mathbb{R}^d} |r_{\log}(\bx,t;\Theta)| p(\bx,t) \dif \bx.
\end{align*}
\end{proof}
\begin{remark}
    In \cref{thm:time-derivative-KL}, we bound the KL divergence between $p(\bx,t)$ and $p_{\Theta}(\bx,t)$ using the residual $r_{\log}(\bx,t;\Theta)$. The residual $r_{\log}(\bx,t;\Theta)$ apparently depends on the modeling capability of $p_{\Theta}(\bx,t)$. Recently some progress on the universal approximation property of invertible mapping has been achieved \cite{Zhu_NN2022,Ishikawa_JMLR2023}. For instance, it has been shown in \cite{Ishikawa_JMLR2023} that the normalizing flow model based on real NVP \cite{dinh2017density} can serve as a universal approximator for an arbitrary PDF in the $L_p$ sense with $p\in[1,\infty)$. Note that our model is a generalization of real NVP. Although tKRnet also has $L_p$-university for PDF approximation, more efforts are needed to establish the convergence with respect to Sobolev norms for the approximation of PDEs, which is beyond the scope of this paper.
\end{remark}

\section{Numerical experiments}
\label{sec:experiments}
To show the effectiveness of our proposed method presented in  \cref{sec:method}, the following four test problems are considered: the double gyre flow problem, the 3-dimensional Kraichnan-Orszag problem, the duffing oscillator, and a 40-variable Lorenz-96 system. 
In our numerical studies, all trainable weights in affine coupling layers (see \eqref{eq:affine-coupling}) are initialized using the Kaiming uniform initialization \cite{he2015delving}, and the biases are set to zero. In each scale-bias layer (see \eqref{eq:scale-bias}), parameters $\ba$ and $\bb$ are initialized as zero, and $\bphi$ is initially set to one. For the nonlinear layer (see \eqref{eq:non-linear-mapping}), the coefficients are initialized as $\hat{m}=32$, $a=50$, and $\psi_i$ and $\varphi_i$ are set to zero for $i=0,\ldots,\hat{m}$ (see \eqref{eq:non-linear-mapping}). The AdamW optimizer \cite{loshchilov2018decoupled} with a learning rate of $0.001$ is used to solve \eqref{eq:opt}, and a cosine annealing learning rate scheduler fine-tunes the learning rate. The training is conducted on an NVIDIA GTX 1080Ti GPU.

In order to assess the accuracy of the tKRnet solution $p_{\Theta}(\bx,t)$ at time $t\in(0,T]$, 
we compare it with the reference solution $p(\bx,t)$ of \eqref{eq:Liouville} as follows, which is computed using the method of characteristics. First, $N_v=10^4$ initial states are sampled from the initial condition $p_0(\bx)$ of \eqref{eq:Liouville}. 
For each initial state, the corresponding solution state of \eqref{eq:uncertain-para-init} at time $t$ is computed using the ordinary differential equation (ODE) solver LSODA in SciPy, and the solution states are denoted by $\{\bx_{\txtval}^{(i)}\}_{i=1}^{N_v}$.     
Then, for each $t\in(0,T]$, the following  relative error is computed  
\begin{equation}
    \label{eq:estimated-rel-abs-err}
    \frac{1}{N_v} \sum_{i=1}^{N_v} \frac{|p(\bx_{\txtval}^{(i)},t)-p_{\Theta}(\bx_{\txtval}^{(i)},t)|}{|p(\bx_{\txtval}^{(i)},t)|},
\end{equation}
and the KL divergence between $p(\bx,t)$ and $p_{\Theta}(\bx,t)$ is estimated as 
\begin{equation}
    \label{eq:estimated-KL}
    D_{KL}(p(\bx,t)||p_{\Theta}(\bx,t))
\approx \frac{1}{N_v}\sum_{i=1}^{N_v} \log\left(\frac{p(\bx_{\txtval}^{(i)},t)}{p_{\Theta}(\bx_{\txtval}^{(i)},t)}\right).
\end{equation} 

\subsection{Double gyre flow}
\label{subsec:Double-gyre-flow}
We start with the nonlinear time-dependent double gyre ﬂow, which has a signiﬁcant effect of nonlinearities for long-time integration  \cite{luchtenburg2014long-time}. The following nonlinear ODE system is considered 
\begin{equation}
\left \{
    \begin{aligned}
        \frac{\dif x_1}{\dif t} =& -\pi A \sin(\pi f(x_1,t)) \cos(\pi x_2)\\
        \frac{\dif x_2}{\dif t} =& \pi A \cos(\pi f(x_1,t))\sin(\pi x_2)\frac{\dif f(x_1, t)}{\dif x_1},
    \end{aligned}
\right .
\end{equation}
where $f(x_1,t) = a(t)x_1^2 + b(t )x_1$, $a(t) = \varepsilon \sin(wt)$, $b(t) = 1-2\varepsilon \sin(wt)$.
In this test problem, the coefficients are set to $A=0.1$, $w=2\pi/10$, and $\varepsilon=0.25$. The time domain is set as $t\in (0,5]$.  
The initial condition $p_0(\bx)$ in \eqref{eq:Liouville} is set to a Gaussian distribution $\mathcal{N}([1,0.5]^\top,0.05^2\mathbb{I})$.  
The tKRnet \eqref{eq:Ttheta} 
has ten time-dependent affine coupling layers, ten scale-bias layers and one nonlinear layer. Each affine coupling layer includes one random Fourier feature layer and two fully connected layers which have thirty two neurons (see \eqref{eq:affine-coupling-structure}). The time domain $(0,5]$ is discretized with step size $\Delta t=0.02$, and the number of spatial collocation points is set to $M=1000$ (see \eqref{eq:residual-time-points}). 
The parameters in \cref{alg:adaptive-sample-flow} are set as $N_r=251000, N_E=100, N_{\adaptive}=6, N_b=251$, and the initial spatial collocation points are generated through the uniform distribution with range $[0,2]\times[0,1]$. 
 
\begin{figure}[htp!]
    \centering
    \subfigure[Reference solution]{\includegraphics[scale=0.26]{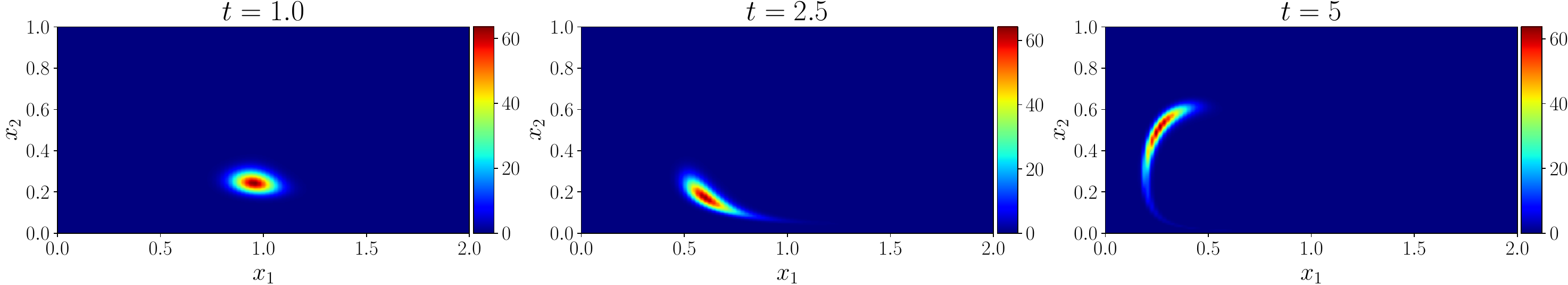}}
    \subfigure[tKRnet solution]{\includegraphics[scale=0.26]{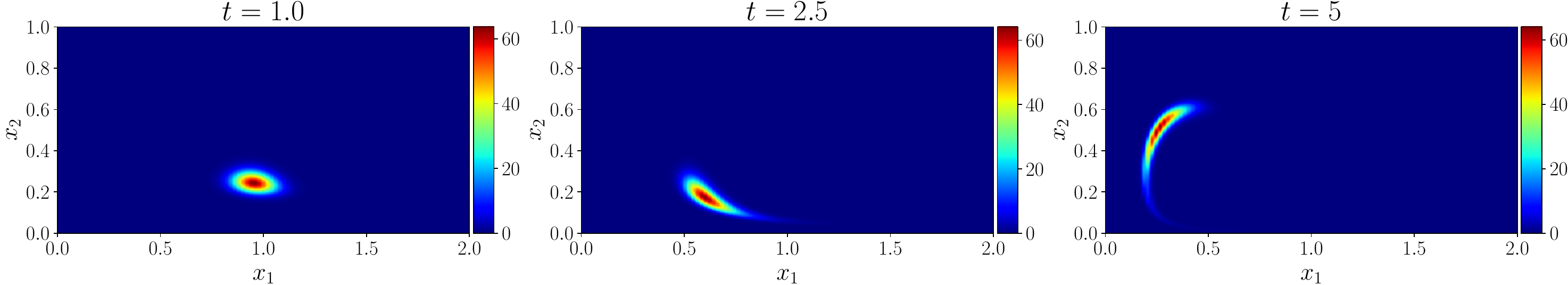}}
    \caption{Double gyre flow problem for $t\in(0,5]$: the reference solution and the tKRnet solution.}
    \label{fig:gyre-flow-solution}
\end{figure}

\begin{figure}[htp!]
    \centering
    \subfigure[Relative absolute error]{\includegraphics[scale=0.3]{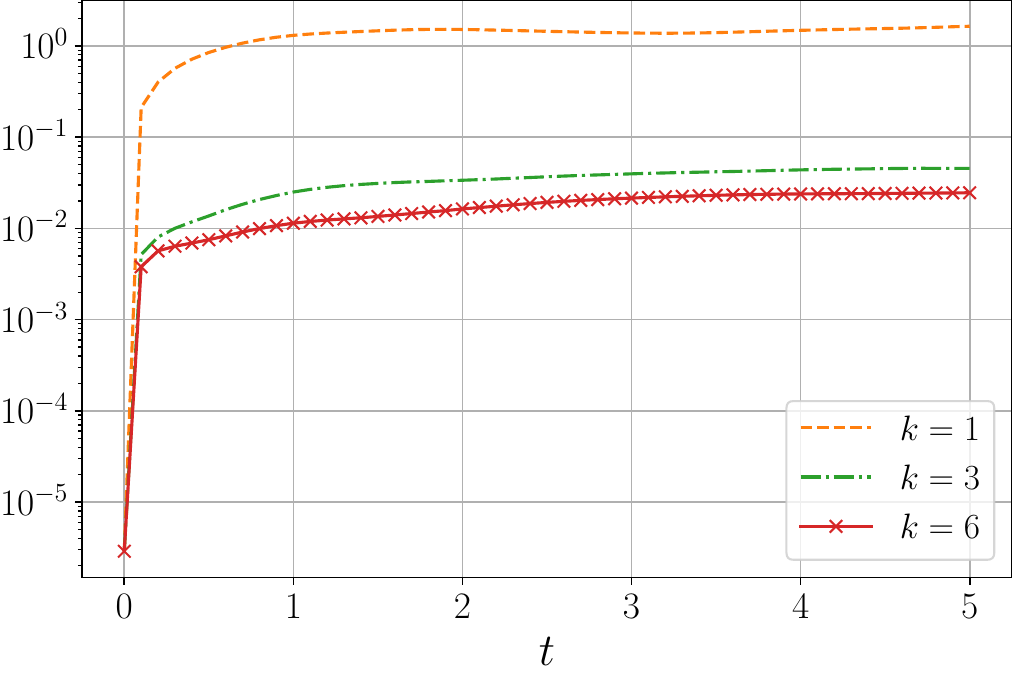}}
    \subfigure[KL divergence]{\includegraphics[scale=0.3]{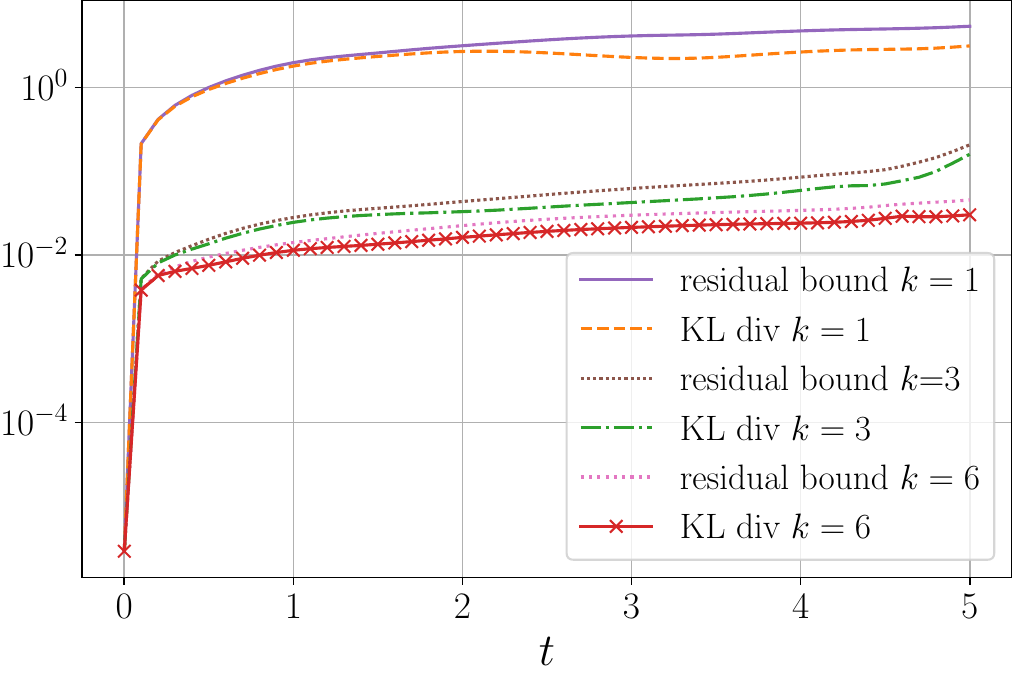}}
    \caption{Double gyre flow problem for $t\in(0,5]$: error and KL divergence values of the tKRnet solutions.}
    \label{fig:gyre-flow-pdf-l2-error}
\end{figure}

{\Cref{fig:gyre-flow-solution}} shows the reference solution obtained using the method of characteristics and our tKRnet solution (obtained using \cref{alg:adaptive-sample-flow}) at three discrete time steps $t=1,2.5,5$.
It can be seen that, 
the tKRnet solution and the reference solution are visually indistinguishable.  
The relative errors (defined in \eqref{eq:estimated-rel-abs-err}) and the values of the KL divergence (defined in \eqref{eq:estimated-KL}) at iteration steps $k=1,3,6$ (see line 4 of \cref{alg:adaptive-sample-flow}) are shown in \cref{fig:gyre-flow-pdf-l2-error}, where it is clear that the errors and the values of the KL divergence significantly reduce as the number of adaptivity iterations increases. In addition, the absolute values of the residual \eqref{eq:log-residual-error} are shown in \cref{fig:gyre-flow-pdf-l2-error}(b). It can be seen that, for each adaptivity iteration step ($k=1,3,6$),  the absolute value of the residual is slightly larger than the value of the KL divergence at each time step, which is consistent with \cref{thm:time-derivative-KL}. 


\begin{figure}[htp!]
    \centering
    \subfigure[Reference solution]{\includegraphics[scale=0.25]{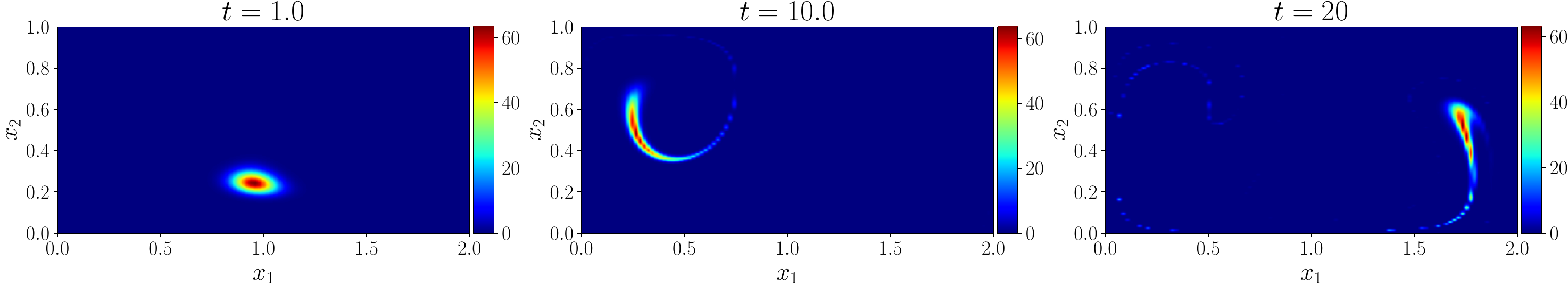}}
    \subfigure[tKRnet solution without temporal decomposition ]{\includegraphics[scale=0.25]{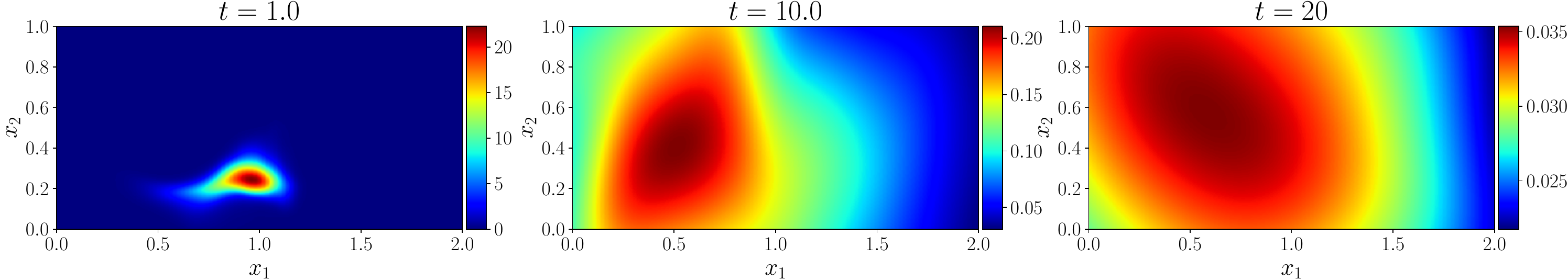}}
    \caption{Double gyre flow problem for $t\in(0,20]$: the reference solution and the  tKRnet solution without temporal decomposition.}
    \label{fig:gyre-flow-long-time-prediction}
\end{figure}

\begin{figure}[!htp]
    \centering
    \subfigure[tKRnet solution with the first choice of temporal decomposition]{\includegraphics[scale=0.26]{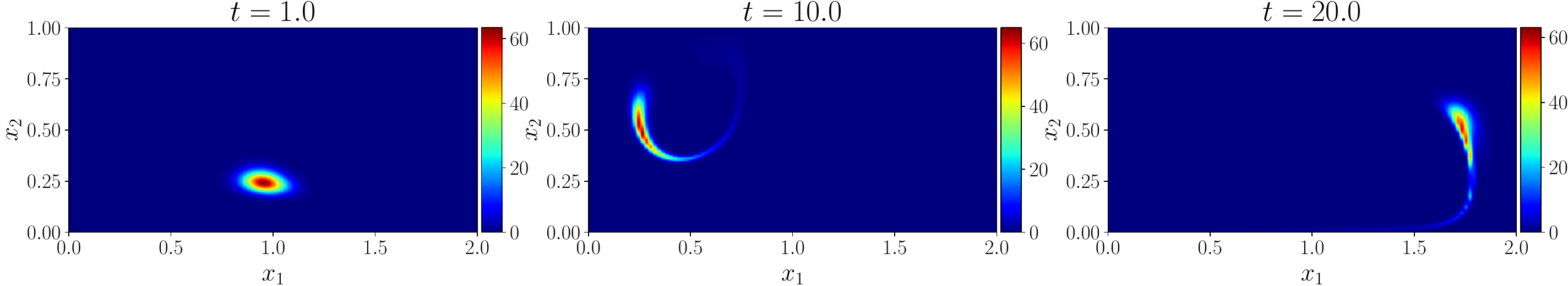}}
    \subfigure[tKRnet solution with the second choice of temporal decomposition]{\includegraphics[scale=0.26]{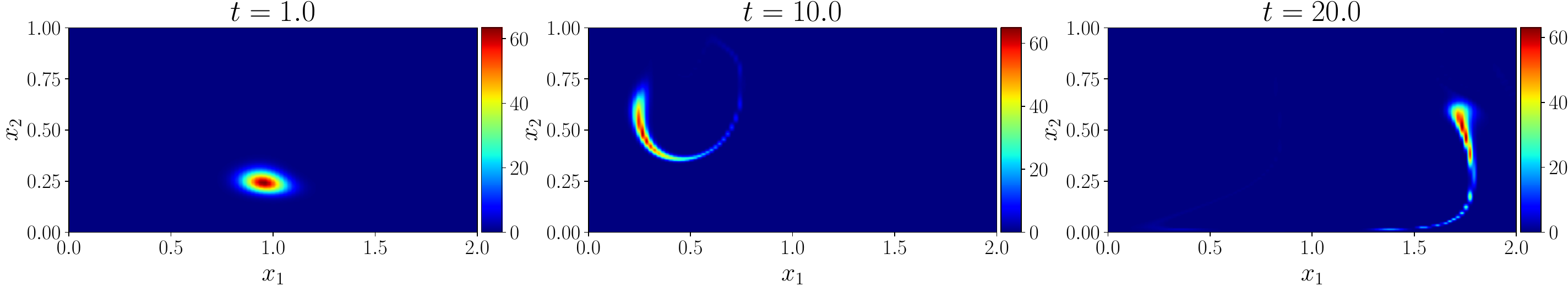}}
    \caption{Double gyre flow problem for $t\in(0,20]$: tKRnet solutions with temporal decomposition.}
    \label{fig:long_time_compared}
\end{figure}

Next, we keep the other settings of this test problem unchanged but extend the time domain to $(0,20]$, which results in a long-time integration problem. The reference solution and our tKRnet solution (directly obtained using \cref{alg:adaptive-sample-flow}) at $t=1,10,20$ are shown in \cref{fig:gyre-flow-long-time-prediction}. It can be seen that, directly applying \cref{alg:adaptive-sample-flow} to this long-time integration problem gives an inaccurate approximation, which is consistent with the challenges addressed in \cite{wang2023long-time}. To resolve this problem, we apply the temporal decomposition method introduced in section \ref{subsec:temporal-dcomposition-long-time}. Here, the interval $(0,20]$ is divided into ten equidistant sub-intervals. Each temporal sub-interval is discretized with time step size $\Delta t=0.02$ ($101$ time steps), and the number of spatial collocation points is set to $M=1000$. So, the total number of collocation points $N_r$ for each sub-interval is $101000=101\times 1000$.
For the first choice in section \ref{subsec:temporal-dcomposition-long-time}, our tKRnet is trained with the loss function \eqref{eq:temporal-decomposition-loss}, and the prior distribution for the tKRnet is set to $\mathcal{N}([1,0.5]^\top,0.05^2\mathbb{I})$. 
For the second choice in section \ref{subsec:temporal-dcomposition-long-time} (see \eqref{eq:prenets-as-prior-approximate-PDF}), 
our tKRnet is trained with the loss function \eqref{eq:loss-func}, where $r_{\log}$ is replaced by $r$ defined \eqref{eq:residual-error} to result in an effective training procedure, 
and the nonlinear layer is not included. \Cref{fig:long_time_compared} shows the results of the two choices for the temporal decomposition. Compared with the reference solution shown in \cref{fig:gyre-flow-long-time-prediction}(a), both choices give efficient tKRnet approximations for this long-time integration problem.


\subsection{Kraichnan-Orszag}
Here we consider the Kraichnan-Orszag problem \cite{WAN2005adaptive}, 
\begin{equation}
\label{eq:Kraichnan-Orszag-system}
\left \{ 
\begin{aligned}
     \frac{\dif x_1}{\dif t}&=x_1 x_3  \\
     \frac{\dif x_2}{\dif t}&=-x_2 x_3 \\
     \frac{\dif x_3}{\dif t}&=-x_1^2 + x_2^2,
\end{aligned}\right. 
\end{equation}
where $x_1$, $x_2$, $x_3$ are state variables. 
The initial condition $p_0(\bx)$ in \eqref{eq:Liouville} is set to a Gaussian distribution $
\mathcal{N}([1,0,0]^\top,0.5^2\mathbb{I})
$, and the time domain in this test problem is set as $t\in(0,3]$.  
The tKRnet \eqref{eq:Ttheta} consists of $\mathT_{[1]}$, $\mathT_{[2]}$ and $L_N$, where each of $\mathT_{[1]},\mathT_{[2]}$ has eight affine coupling layers and eight scale-bias layers. Each affine coupling layer includes one random Fourier layer and three fully connected layers with thirty-two neurons (see \eqref{eq:affine-coupling-structure}). The time domain is discretized with time step size $\Delta t=0.01$, and the number of spatial collocation points is set to $M=4000$ (see \eqref{eq:residual-time-points}).
The settings in  \cref{alg:adaptive-sample-flow} are set as $N_r=1204000, N_E=50, N_{\adaptive}=10, N_b=1204$, and initial spatial collocation points are dawn from 
 the uniform distribution with range $[-5,5]^3$.

\begin{figure}[htp!]
   \centering
    \subfigure[Reference solution]{
      \includegraphics[scale=0.3]{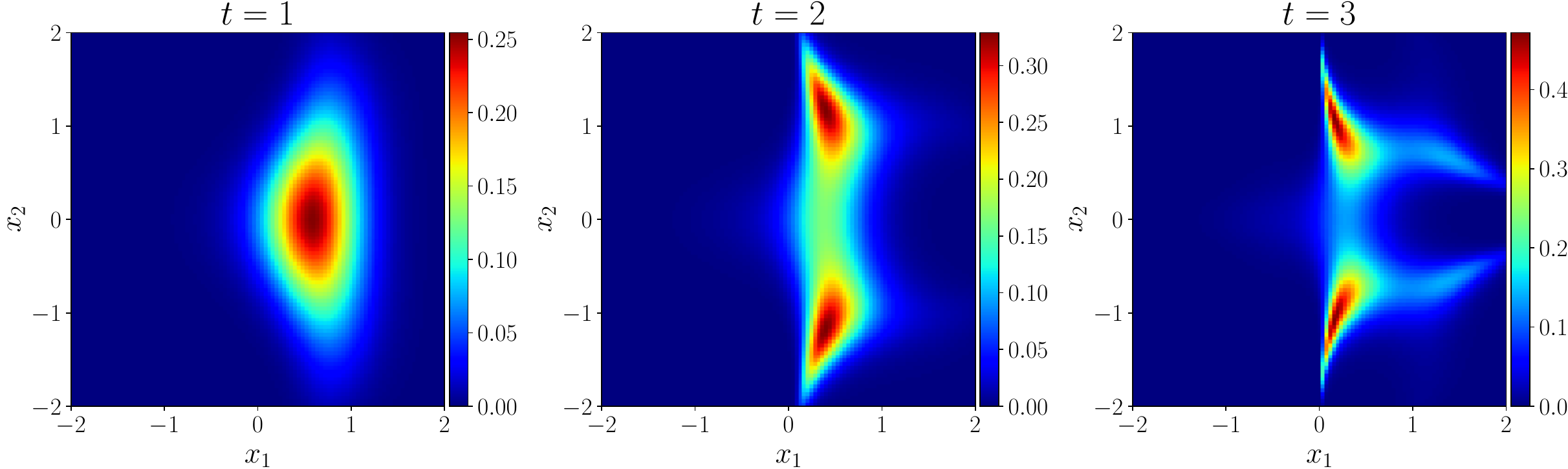}
      \label{subfig:Kraichnan-Orszag-ref-sol}
      }
    \subfigure[tKRnet solution]{
      \includegraphics[scale=0.3]{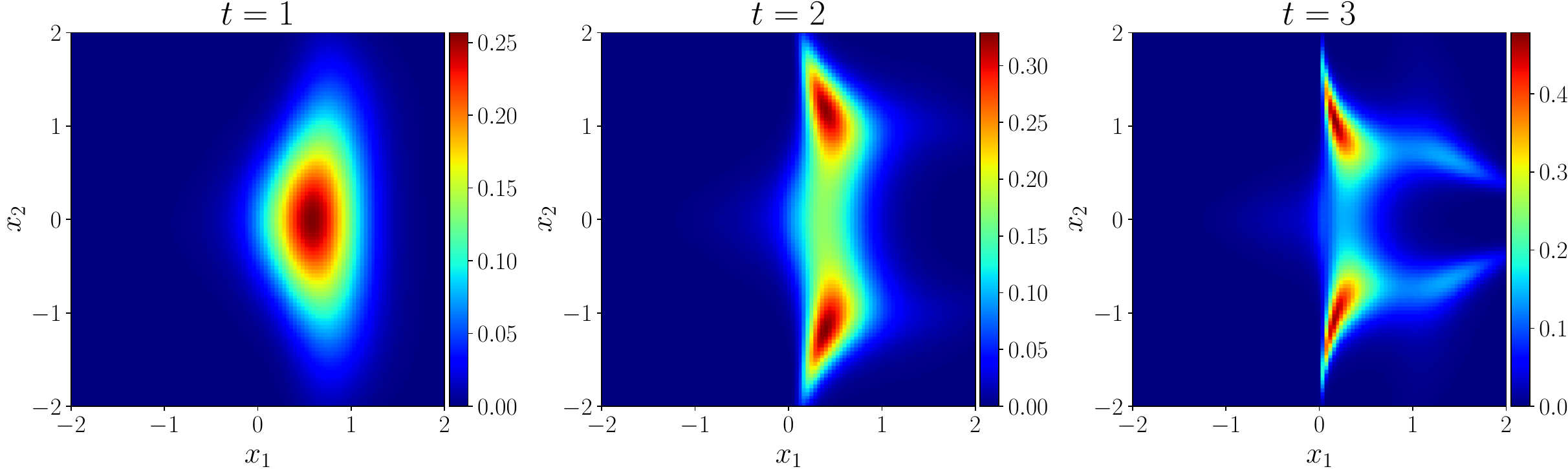}
      \label{subfig:Kraichnan-Orszag-approximate-sol}
      } 
\caption{Kraichnan-Orszag problem: The reference solution and the tKRnet solution.}
   \label{fig:Kraichnan-Orszag-PDF}
\end{figure}

\begin{figure}[htp!]
    \centering
    \subfigure[Relative absolute error]{\includegraphics[scale=0.3]{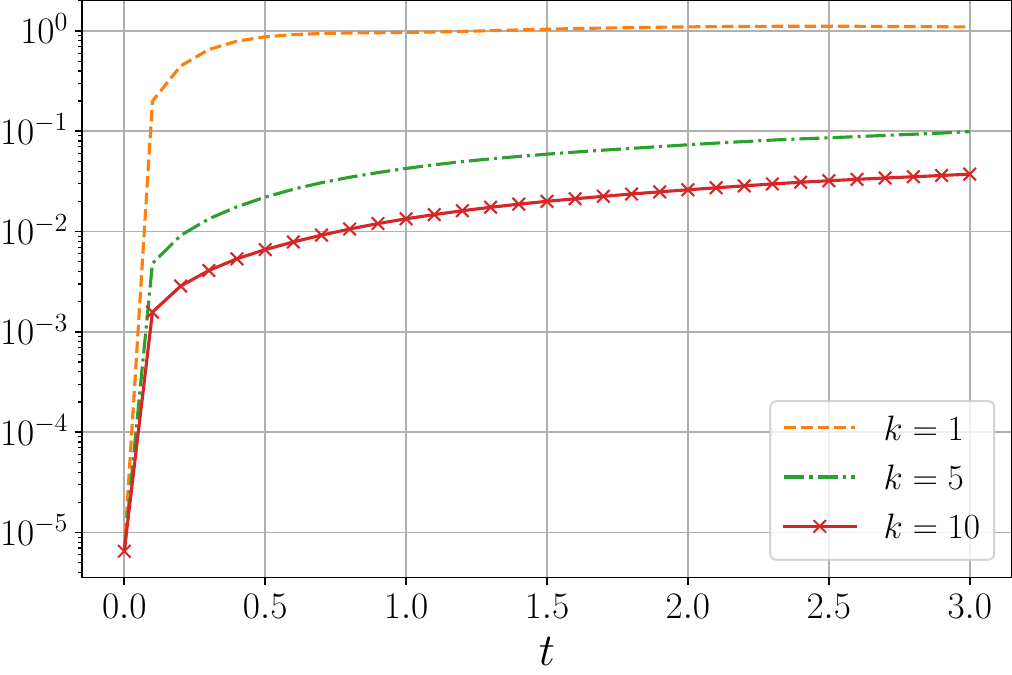}}
    \subfigure[KL divergence]{\includegraphics[scale=0.3]{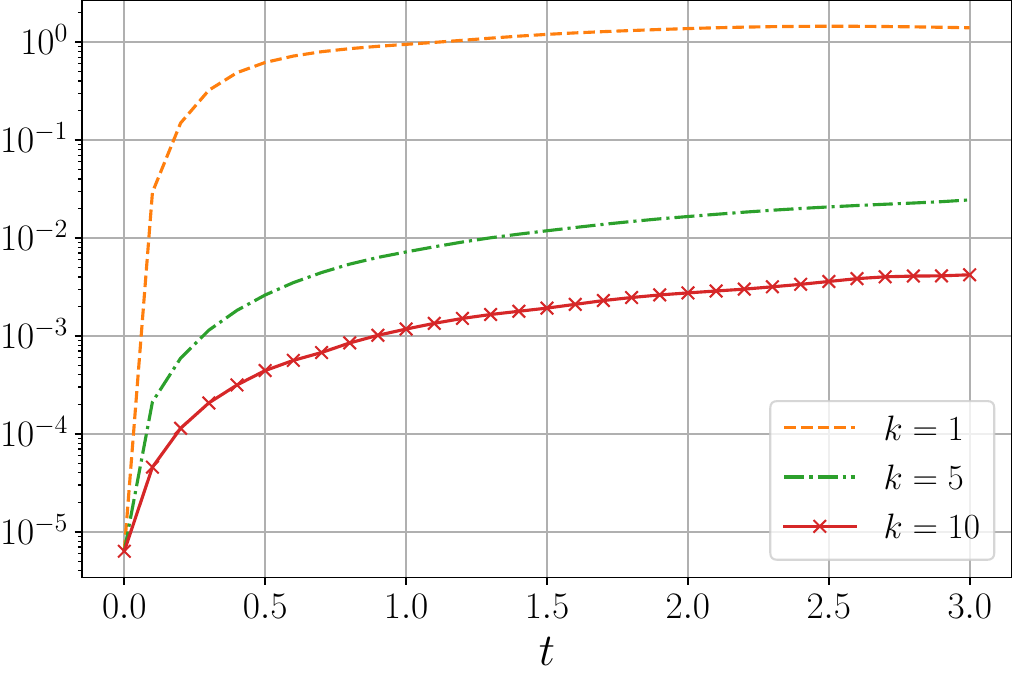}}
    \caption{Kraichnan-Orszag problem: error and KL divergence values of the tKRnet solutions.}
    \label{fig:Kraichnan-Orszag-pdf-l2-error}
\end{figure}

\Cref{fig:Kraichnan-Orszag-PDF} shows the tKRnet solution and the reference solution  at $t=1,2,3$, where it can be seen that the tKRnet solution 
and the reference solution are visually indistinguishable. 
The values of the relative error \eqref{eq:estimated-rel-abs-err} and the KL divergence \eqref{eq:estimated-KL} at three adaptivity iterations $k=1,5,10$ (see line 4 of \cref{alg:adaptive-sample-flow}) are illustrated in \cref{fig:Kraichnan-Orszag-pdf-l2-error}.  
It is clear that the errors and the values of the KL divergence decrease as the number of adaptivity iteration steps increases.  


\subsection{Forced Duffing oscillator}
The forced Duffing oscillator system is defined as:
\begin{equation}
    \label{eq:duffing oscillator}
    \left \{
    \begin{aligned}
        \frac{\dif y_1}{\dif t} &=y_2 \\
    \frac{\dif y_2}{\dif t} &=-\delta y_2 - y_1 (\alpha+\beta y_1^2) + \gamma \cos(\omega t),
    \end{aligned}
    \right.
\end{equation}
where $y_1$ and $y_2$ represent the state variables and $\delta$, $\alpha$, $\beta$, $\gamma$, and $\omega$ represent uncertain parameters. The time domain is set to $(0,2]$. The initial distribution of the state variables $p_{\by}(\by,0)$ and the distribution of the uncertain parameters $p_{\bxi}(\bxi)$ in \eqref{eq:Liouville} are set as a Gaussian distribution $\mathcal{N}([0,0]^\top,\mathbb{I})$ and a Gaussian distribution $\mathcal{N}([0.5,-1,1,0.5,1]^\top,0.25^2\mathbb{I})$ respectively. Letting $\bx=[\by,\bxi]^\top$, the initial condition 
in \eqref{eq:Liouville} is constructed as $p_0(\bx)=p_{\by}(\by,0)p_{\bxi}(\bxi)$.  
The tKRnet \eqref{eq:Ttheta} consists of $\mathT_{[1]},\mathT_{[2]},\mathT_{[3]}$ and $L_N$, where each of $\mathT_{[1]},\mathT_{[2]}$ and $\mathT_{[3]}$ has four affine coupling layers and four scale-bias layers. Each affine coupling layer has one random Fourier layer and two fully connected layers with thirty two neurons (see \eqref{eq:affine-coupling-structure}). The coefficients in \cref{alg:adaptive-sample-flow} are set as $N_r=804000, N_E=100, N_{\adaptive}=6, N_b=804$. The time domain is discretized with time step size $\Delta t=0.01$, and the number of spatial collocation points is $M=4000$ (see \eqref{eq:residual-time-points}). Initial spatial collocation points are sampled from the uniform distribution with range $[-5,5]^7$. 

\begin{figure}[htp!]
    \centering
    \subfigure[Reference solution]{ \includegraphics[scale=0.3]{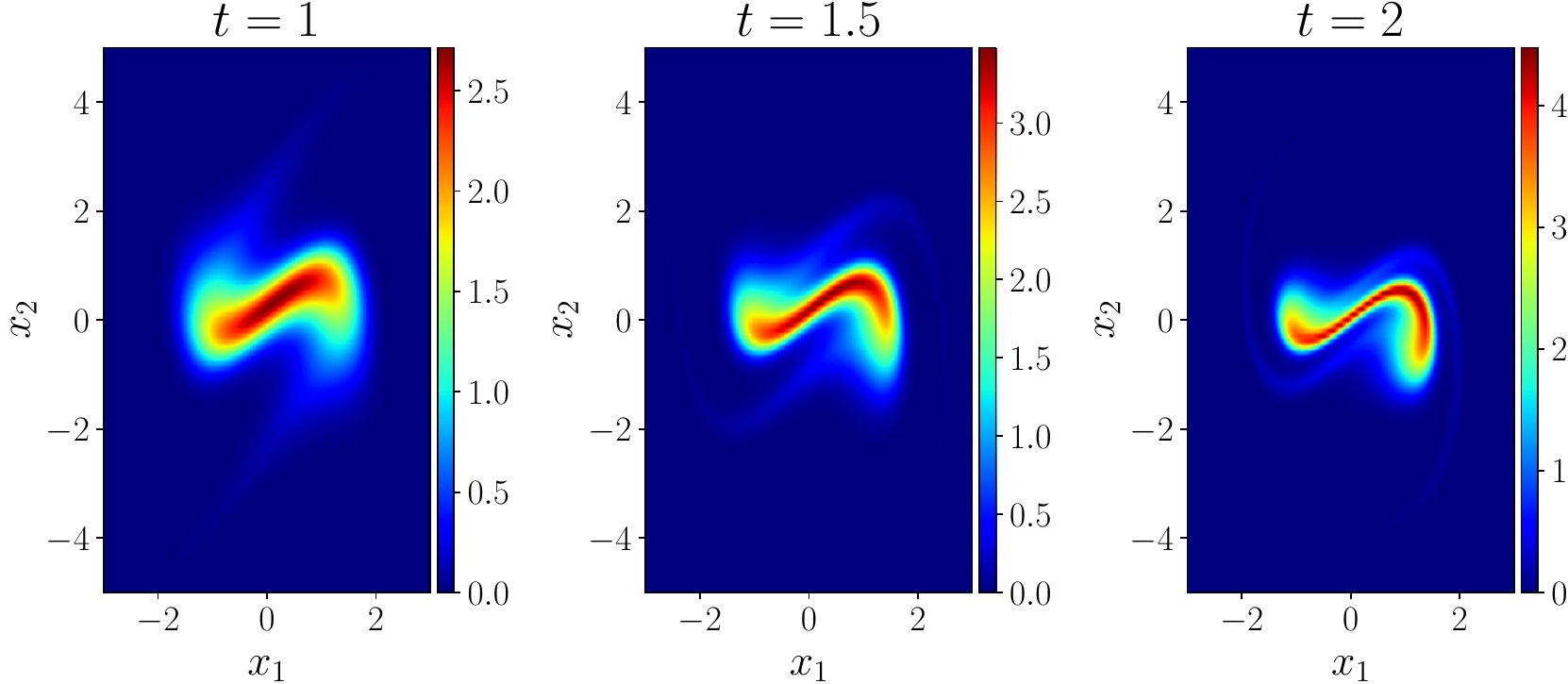}}
   \subfigure[tKRnet solution]{ \includegraphics[scale=0.3]{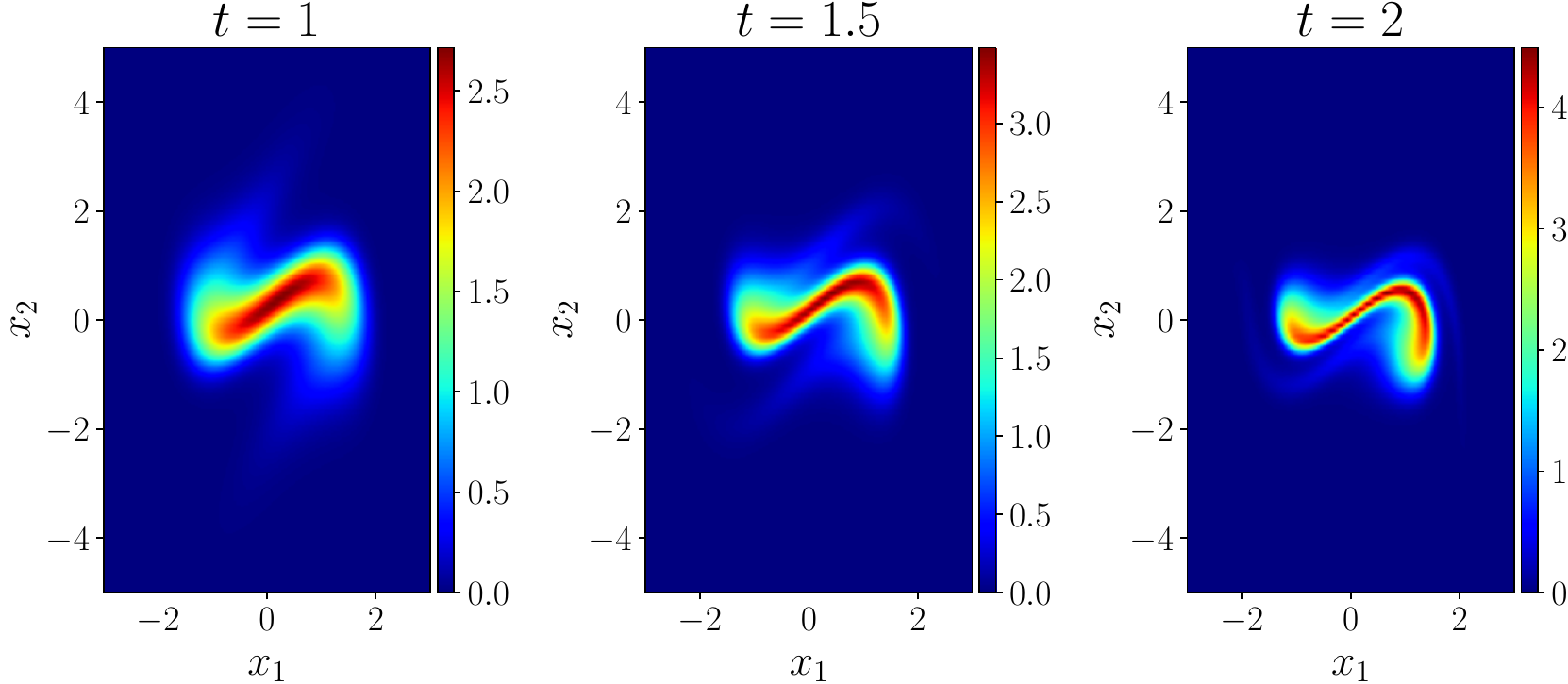}
   \label{subfig:duffing-oscillator-approx-sol}
   }
    \caption{Duffing oscillator problem: the reference solution and the tKRnet solution.}
    \label{fig:duffing-oscillator-ref-sol-approximate-sol}
\end{figure}

\begin{figure}[htp!]
    \centering
     \subfigure[Relative absolute error]{\includegraphics[scale=0.3]{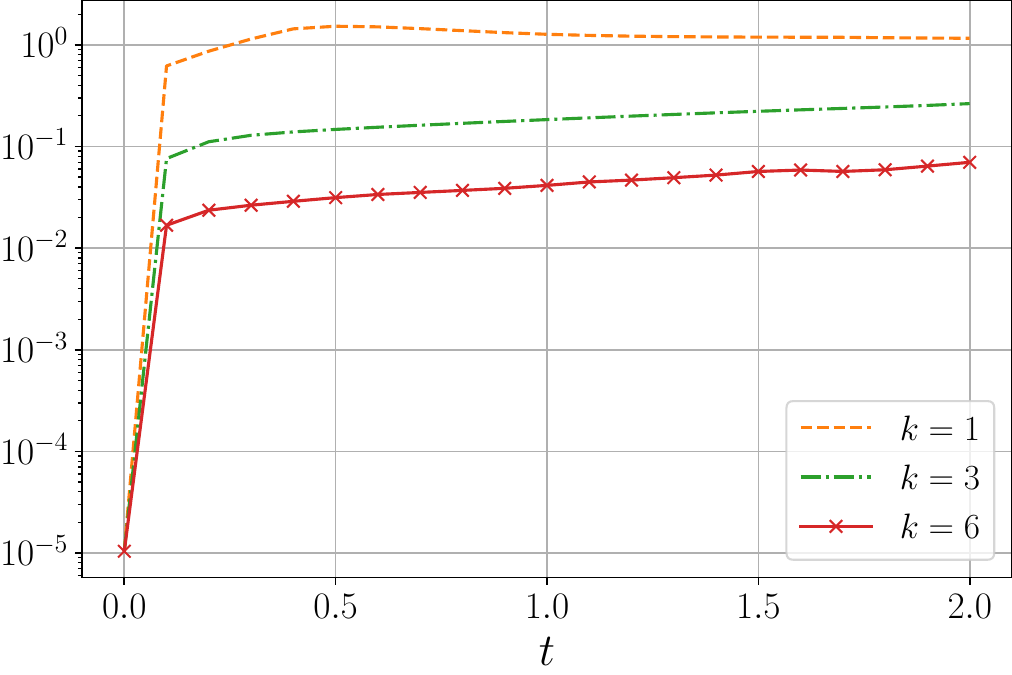}}
    \subfigure[KL divergence]{\includegraphics[scale=0.3]{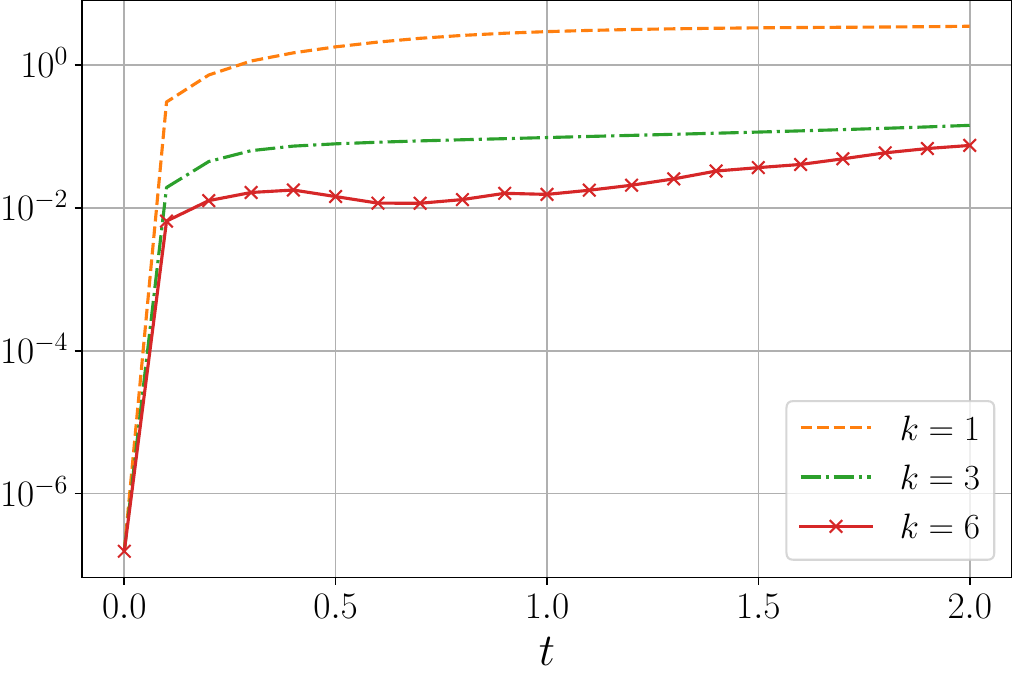}}
    \caption{Duffing oscillator problem: 
    error and KL divergence values of the tKRnet solutions.}
    \label{fig:duffing-oscillator-KL-div-mae}
\end{figure}

\Cref{fig:duffing-oscillator-ref-sol-approximate-sol} shows the tKRnet solution and the reference solution at $t=1,1.5,2$. 
From \cref{fig:duffing-oscillator-ref-sol-approximate-sol}, it can be seen that the tKRnet solution and the reference solution are visually indistinguishable. 
The relative absolute errors \eqref{eq:estimated-rel-abs-err} and the values of KL divergence \eqref{eq:estimated-KL} at three adaptivity iterations $k=1,3,6$ (see line 4 of \cref{alg:adaptive-sample-flow}) are illustrated in \Cref{fig:duffing-oscillator-KL-div-mae}.  
It is clear that the errors and the values of KL divergence decrease as the number of adaptivity iterations increases. 
 
\subsection{Lorenz-96 system}
In this test problem, the Lorenz-96 system is considered, which is a model used in numerical weather forecasting \cite{Karimi2010extensive}. 
The general form of the Lorenz-96 system is defined as  
\begin{equation}
    \label{eq:lorenz-96}
    \frac{\dif x_i}{\dif t} = (x_{i+1}-x_{i-2})x_{i-1}-x_i+F, \quad i=1,\ldots,d,
\end{equation}
where $x_i$ (for $i=-1,\ldots, d+1$) are the state variables, and $F$ represents a constant force. 
In this test problem, it is assumed that $d\geq 4$, $x_{-1}=x_{d-1}$, $x_0=x_d$ and $x_{d+1}=x_1$. We set  $d=40$, $F=1$, and $t\in(0,1]$. The initial condition $p_0(\bx)$ in \eqref{eq:Liouville} is set to the joint Gaussian distribution 
\begin{align*}
    p_0([x_1,\ldots,x_{40}]^\top)=\left(\frac{25}{2\pi}\right)^{20} \prod_{i=1}^{40}\exp\left(-\frac{25}{2}(x_i-(0.5-|\frac{i}{40}-0.5|))^2\right).
\end{align*}
The tKRnet \eqref{eq:Ttheta} for this problem has a sequence of transformations $\mathT_{[1]},\ldots,\mathT_{[5]}$ and one nonlinear layer $L_N$, where each $\mathT_{[1]},\ldots,\mathT_{[5]}$ includes four affine coupling layers and four scale-bias layers. Each affine coupling layer has one random Fourier layer and two fully connected layers with 128 neurons. The time domain is discretized with time step size $\Delta t=0.01$, and the number of spatial collocation points is set to $M=2000$ (see \eqref{eq:residual-time-points}). 
Parameters in \cref{alg:adaptive-sample-flow} are set as 
$N_r=202000, N_E=50, N_{\adaptive}=10, N_b=202$, and initial spatial collocation points are generated using the uniform distribution with range $[-5,5]^{40}$. 

\begin{figure}[htp!]
    \centering
    \subfigure[Reference solution]{\includegraphics[scale=0.34]{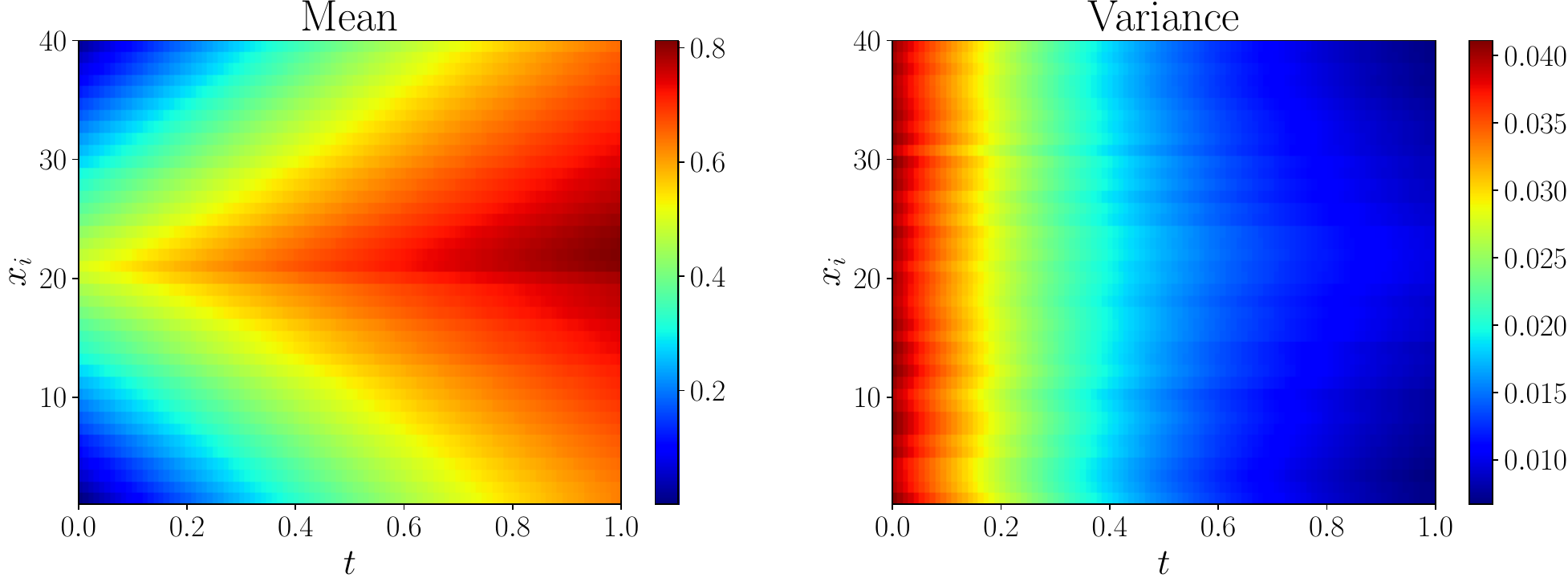}}
    \subfigure[tKRnet solution]{\includegraphics[scale=0.34]{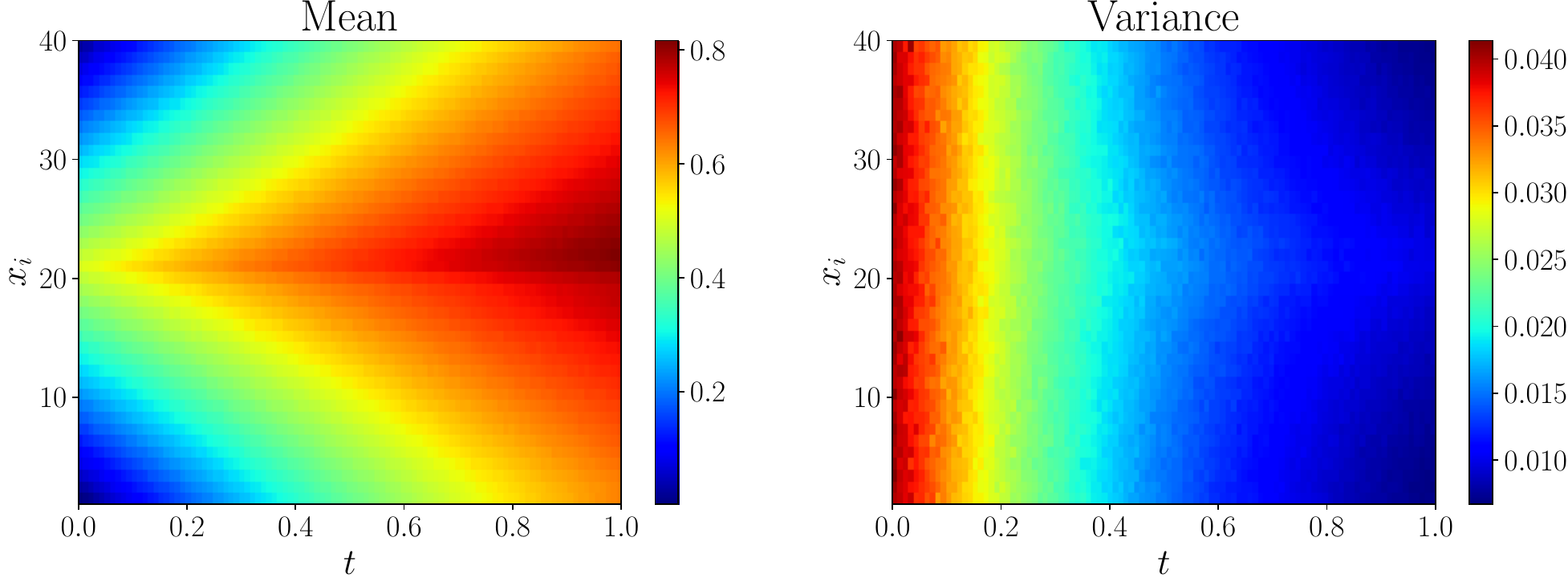}}
    \caption{Lorenz-96 problem: mean and variance of reference and tKRnet solutions.}
    \label{fig:lorenz-96-mean-var}
\end{figure}

\begin{figure}[htp!]
    \centering
    \includegraphics[scale=0.34]{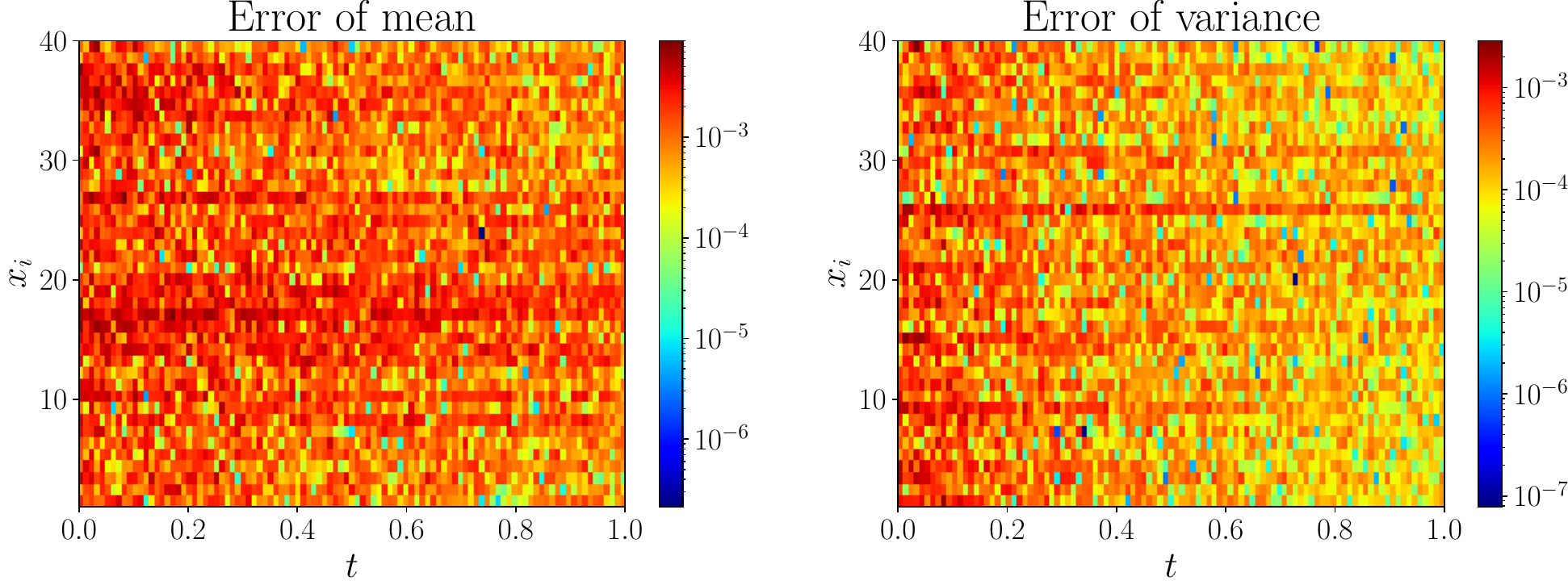}
    \caption{Lorenz-96 problem: the absolute error to mean and variance in log scale.}
    \label{fig:lorenz-96-mean-var-error}
\end{figure}

For this high-dimensional problem, 
the mean and the variance estimates of the reference solution and the tKRnet solution are compared. For the reference solution, $N_{v}=10^4$ initial states are sampled from $p_0(\bx)$, and the states $\{\bx_{\txtval}^{(i)}\}_{i=1}^{N_{v}}$ at time $t\in(0,T=1]$ are obtained by solving \eqref{eq:uncertain-para-init} using the LSODA solver. The mean and the variance estimates are computed as 
\begin{eqnarray}
    \widehat{\mathbb{E}}_{p(\bx,t)}[\bx;t] &:=& \frac{1}{N_{v}}\sum_{i=1}^{N_{v}}\bx_{\txtval}^{(i)},\label{eq_mean}\\
    \widehat{\Var}_{p(\bx,t)}[\bx;t] &:=& \frac{N_{v}}{N_{v}-1}\left( \frac{1}{N_{v}}\sum_{i=1}^{N_{v}}(\bx_{\txtval}^{(i)})^2-\left(\frac{1}{N_{v}}\sum_{i=1}^{N_{v}}\bx_{\txtval}^{(i)}\right)^2\right). \label{eq_var}
\end{eqnarray} 
For tKRnet approximation solution, for a given time $t\in(0,1]$, samples of the states are generated by $p_{\Theta}(\bx,t)$, and the mean and variance estimates are obtained by putting the samples into \eqref{eq_mean} and \eqref{eq_var}, which are denoted by $\widehat{\mathbb{E}}_{p_{\Theta}(\bx,t)}$ and $\widehat{\Var}_{p_{\Theta}(\bx,t)}$ respectively.  \Cref{fig:lorenz-96-mean-var} shows the reference mean and variance estimates of the reference solution and the tKRnet solution, where it can be seen that the results of the reference solution and those of the tKRnet solution are very close. 
Next, at time $t\in(0,1]$, the errors in the mean and variance estimates are computed as 
\begin{eqnarray*}
\left|\widehat{\mathbb{E}}_{p(\bx,t)}[\bx;t]-\widehat{\mathbb{E}}_{p_{\Theta}(\bx,t)}[\bx;t]\right|, \\
\left|\widehat{\Var}_{p(\bx,t)}[\bx;t]-\widehat{\Var}_{p_{\Theta}(\bx,t)}[\bx;t]\right|.
\end{eqnarray*}
\Cref{fig:lorenz-96-mean-var-error} shows the errors, where it is clear that the errors are small---the maximum of the errors in the mean estimate is around $7.420\times 10^{-3}$ and that in the variance estimate is around $ 2.657 \times 10^{-3}$.

\section{Conclusions}
\label{sec:conclusions}
The uncertainty quantification of stochastic dynamical systems can be addressed by computing the time-dependent PDF of the states. However, the states of the system can be high-dimensional and the support of the PDF may be unbounded. To address these issues, we have proposed a physics-informed adaptive density approximation method based on tKRnets to approximate the Liouville equations. The tKRnet provides an explicit family of PDFs via the change of variable rule with a trainable time-dependent invertible transformation. The initial PDF of the Liouville equation can be encoded in the tKRnet through the prior distribution. Adaptive sampling plays a crucial role in achieving an accurate PDF approximation, where the training set needs to be updated according to the localized information in the solution. 
By coupling the adaptive sampling with an efficient temporal decomposition, the long-time integration can be effectively improved. Numerical results have demonstrated the efficiency of our algorithm for high-dimensional stochastic dynamical systems. In this work, we use uniform grids to discretize the time interval without paying much attention to the causality in the time direction. Adaptivity can be introduced into temporal discretization for further refinement. This issue is being investigated and will be reported elsewhere. 

\appendix
\section{Additional training results with ODE residual}
The PDEs \cref{eq:Liouville} and \cref{eq:log-Liouville} can be solved using the method of characteristics \cite{morton2005numerical}, where the characteristic lines evolve along the solution of the stochastic ODE system \cref{eq:uncertain-para-init}. Therefore, an alternative approach to learn the solution of \cref{eq:Liouville} (or \cref{eq:log-Liouville}) is to construct a time-dependent invertible mapping $\mathT$ (a deep neural network) to approximate the flow map for \cref{eq:uncertain-para-init}. 
We define $\bz = \mathT(\bx,t;\Theta)$ and its inverse $\bx = \mathT^{-1}(\bz,t;\Theta)$, where $\bz\sim p_0$, $\bx$ is the state variable in \cref{eq:uncertain-para-init} and $\Theta$ is the parameters of the neural network $\mathT$. To ensure $\mathT^{-1}$ learns the solution of \eqref{eq:uncertain-para-init}, the following residual is defined
\begin{equation*}
    \ell_{ode}(\bz, t;\Theta) = \left\|\frac{\partial \mathT^{-1}(\bz, t;\Theta)}{\partial t} - \mathf(\mathT^{-1}(\bz, t;\Theta), t)\right\|_2^2,\textrm{ where } \bz=\mathT(\bx, t;\Theta);
\end{equation*}
the corresponding loss function is given by 
\begin{equation}
    \label{eq:ode-loss}
    \begin{aligned}
        \mathcal{L}_{ode}(\bz, t;\Theta) = \frac{1}{N_r}\sum_{i=1}^{N_r} \ell_{ode}(\bz_{\text{res}}^{(i)}, \rest^{(i)};\Theta),\, \textrm{ where } \bz_{\text{res}}^{(i)}=\mathT(\resx^{(i)}, \rest^{(i)};\Theta),
    \end{aligned}
\end{equation}
where $\{\resx^{(i)}\}_{i=1}^{N_r}$ are spatial collocation points  
and $\{\rest^{(i)}\}_{i=1}^{N_r}$ are temporal collocation points (see \eqref{eq:residual-time-points}). Then, the deep neural network $\mathT$ can be trained using \cref{alg:adaptive-sample-flow} with the loss function \eqref{eq:ode-loss}, 
and the ODE based approximation is constructed as $p_0(\mathT(\bx,t;\Theta))|\det \nabla_{\bx}\mathT(\bx,t;\Theta)|$ using the trained neural network $\mathT(\bx,t;\Theta)$, while the PDE based approximation is the approximate PDF by minimizing \eqref{eq:loss-func}. 

\begin{figure}[htp!]
    \centering
    \includegraphics[scale=0.26]{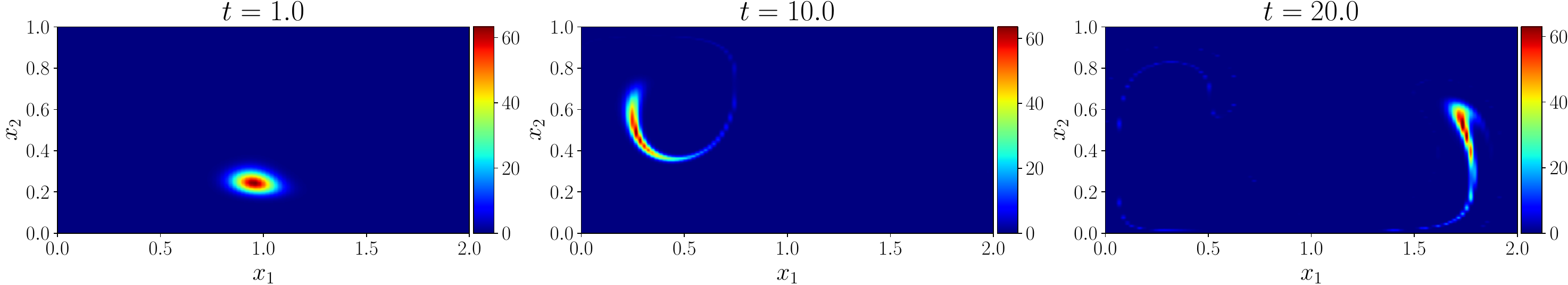}
    \caption{ODE based approximation.}
    \label{fig:ode-loss-long-time-pred}
\end{figure}

\begin{figure}[htp!]
    \centering
    \includegraphics[scale=0.32]{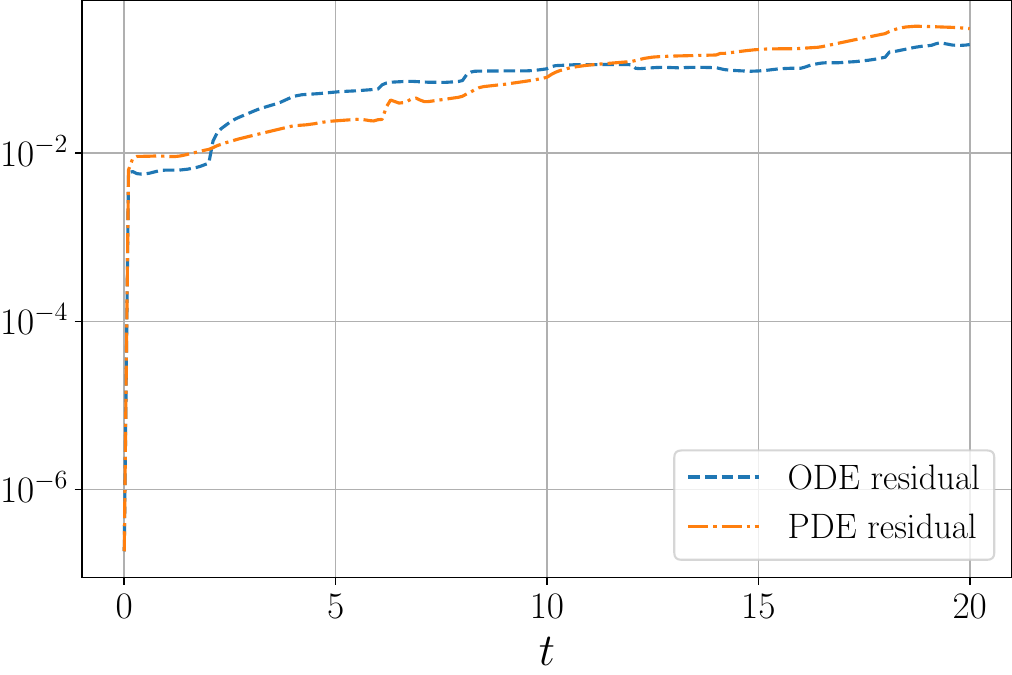}
    \caption{Relative errors of ODE based approximation and PDE based approximation.}
    \label{fig:ODE-PDE-residual-comparison}
\end{figure}

For the long-time integration problem considered in section \ref{subsec:Double-gyre-flow}, the ODE based approximation and the PDE based approximation are compared as follows, where the second choice for temporal decomposition (introduced in section \ref{subsec:temporal-dcomposition-long-time}) is applied to both approximations. All settings are the same as those in section \ref{subsec:Double-gyre-flow} for $t\in(0,20]$. \cref{fig:ode-loss-long-time-pred} shows the ODE based approximation, which approximates the reference solution (\cref{fig:gyre-flow-long-time-prediction}(a)) well. 
\Cref{fig:ODE-PDE-residual-comparison} shows relative errors (defined in \eqref{eq:estimated-rel-abs-err}) of ODE and PDE approximations, where it is clear that the relative errors of both approximations are comparable. However, as the training procedure with the loss function \eqref{eq:ode-loss} requires computing ${\partial \mathT^{-1}(\bz, t;\Theta)}/{\partial t}$ with backpropagation, the cost for training the ODE based approximation is significantly larger than that for training the PDE based approximation, especially when the state is high-dimensional. 


\bibliographystyle{siamplain}
\bibliography{references}
\end{document}